\documentclass[11pt,a4paper]{article}
\usepackage[T1]{fontenc}
\usepackage[utf8]{inputenc}
\usepackage{lmodern}
\usepackage[a4paper,margin=27mm]{geometry}
\usepackage{amsmath,amssymb,amsthm,mathtools}
\usepackage{microtype}
\usepackage{enumitem}
\usepackage{booktabs,array}
\usepackage[hidelinks,unicode]{hyperref}
\usepackage{bookmark}
\newcommand{\Z}{\mathbb{Z}}
\newcommand{\R}{\mathbb{R}}

\newtheorem{theorem}{Theorem}[section]
\newtheorem{lemma}[theorem]{Lemma}
\newtheorem{proposition}[theorem]{Proposition}
\newtheorem{corollary}[theorem]{Corollary}
\theoremstyle{definition}
\newtheorem{definition}[theorem]{Definition}
\newtheorem{example}[theorem]{Example}

\newtheorem{mainresult}{Main Result}

\numberwithin{equation}{section}
\setlist[enumerate]{itemsep=3pt,topsep=5pt}
\allowdisplaybreaks[2]
\title{Universal optimality on finite Lee tori}
\author{Tao Zhang\thanks{Email: \texttt{zhant220@163.com}.}\\
    \small Institute of Mathematics and Interdisciplinary Sciences,\\[-1mm]
    \small Xidian University, Xi'an 710126, China}
\date{}
\begin{document}
    \maketitle

    \begin{abstract}
        We study subsets of finite Lee tori that minimize pair energy among sets of the same cardinality for every discretely completely monotone potential.
        At half occupancy on a product of even cycles, the two checkerboards minimize both total pair energy and maximum individual energy. They are the only minimizers for every nonconstant classical completely monotone potential on $(0,\infty)$, proving the checkerboard conjecture of Bouman, Draisma, and van Leeuwaarden for this class.
        At the independence number of the distance one graph on $\Z_q^n$, universally optimal sets exist exactly when $q$ is even, $q=3$, $n=1$, or $(q,n)=(5,2)$. In these cases they are precisely the maximum independent sets and form a single Lee isometry class.
        At density $1/5$ in $\Z_5^n$, $n\ge2$, universal optimality is characterized by a weighted count of pairs at distances one and two.
        We prove these results using integer box occupancy bounds that strengthen Fourier linear programming and, in the quinary case, Fourier certificates with a common zero set.

        
        
    \end{abstract}

    \section{Introduction}

    A pair potential assigns a cost $f(r)$ to two points at distance $r$. For an $M$-point subset $C$ of a metric space $(X,d)$, the mean individual energy is
    \[
    E_f(C)=\frac{1}{M}\sum_{x\in C}\sum_{y\in C\setminus\{x\}}f(d(x,y)).
    \]
    As the exponent $s$ increases in the inverse power law $f(r)=r^{-s}$, the shortest distances carry an increasing share of the energy. A basic question is whether a single configuration can minimize energy for every exponent, or more generally for an entire class of repulsive potentials.

    Cohn and Kumar~\cite{ref9} introduced the notion of a \emph{universally optimal} configuration on a sphere, it minimizes energy among configurations of the same cardinality for every completely monotone function of squared Euclidean distance. Using a linear programming bound for spherical energy~\cite{ref25,ref9}, they proved, for example, that the regular icosahedron has this property among twelve-point subsets of the unit sphere in $\R^3$. The Euclidean linear programming method for sphere packing~\cite{ref8} also underlies the proof of universal optimality for the $E_8$ and Leech lattices at fixed density~\cite{ref10}.

    For codes in Hamming space, Cohn and Zhao~\cite{ref12} formulated complete monotonicity using finite differences in place of derivatives. On $\{1,\ldots,D\}$, a potential is \emph{completely monotone} if $(-1)^j\Delta^j f(r)\ge0$ whenever the difference is defined, where $\Delta f(r)=f(r+1)-f(r)$. Inverse powers $r^{-s}$, $s>0$, and exponentials $c^r$, $0<c<1$, are examples. The binomial potentials $b_j(r)=\binom{D-r}{j}$, $0\le j<D$, generate this cone. Universal optimality is therefore equivalent to minimizing the corresponding binomial moments of the distance distribution. Bounds on these moments~\cite{ref1}, together with Delsarte's linear programming method~\cite{ref13,ref14}, yield many examples in Hamming space~\cite{ref12,ref6}.

    We consider a finite abelian group $G=\prod_{i=1}^n\Z_{q_i}$ with specified cyclic factors of orders $q_i\ge2$. Write $V=|G|$ and $D=\sum_i\lfloor q_i/2\rfloor$. For $0\le a<q$, put $\ell_q(a)=\min(a,q-a)$. The Lee distance~\cite{ref20} is the graph distance on the torus,
    \[
    d_L(x,y)=\sum_{i=1}^n\ell_{q_i}(x_i-y_i).
    \]
    Lee and Hamming distances agree on $\Z_2^n$ and $\Z_3^n$, and the Gray map identifies $\Z_4^n$ with a binary Hamming space~\cite{ref17}. The distinction between the two metrics already appears in $\Z_5^2$, where two lines that are equivalent as Hamming codes can have different Lee energies. The distances from any point of $y=-x$ to the other four points are $2,2,4,4$, whereas on $y=2x$ all four distances are $3$. Strict convexity thus favors the second line (Example~\ref{ex:index-p-failure}). The Lee association scheme~\cite{ref2,ref23} records the coordinate lengths of a displacement up to permutation, a distinction also used in linear and semidefinite bounds for Lee codes~\cite{ref3,ref22}.

    Bouman, Draisma, and van Leeuwaarden~\cite{ref5} studied energy minimization under the Lee metric in a model of wireless interference. Their objective is the greatest interference at an occupied site, namely the \emph{maximum individual energy}
    \[
    E_f^{\max}(C)=\max_{x\in C}\sum_{y\in C\setminus\{x\}}f(d_L(x,y)).
    \]
    When all side lengths are even, the torus splits into two checkerboards. Conjecture~1.1 of~\cite{ref5} asks whether they are the only minimizers of $E_f^{\max}$ at half occupancy for a suitable class of repulsive potentials. For smooth strictly completely monotone potentials, the authors proved this when each side length is two or a multiple of four~\cite[Theorem~5.3]{ref5}. Brosch and de Klerk obtained semidefinite bounds for specific grids~\cite{ref7}. Our first result treats arbitrary even side lengths and gives an explicit condition for uniqueness.

    \begin{mainresult}[Theorem~\ref{thm:checkerboard}]\label{main:checkerboard}
        If every $q_i$ is even, the two checkerboards $C_\varepsilon=\{x\in G:\sum_i x_i\equiv\varepsilon\pmod2\}$, $\varepsilon\in\{0,1\}$, minimize both $E_f$ and $E_f^{\max}$ among all $V/2$-point sets for every completely monotone $f$. They are the only minimizers of either energy if $D\ge2$ and $(-1)^{D-1}\Delta^{D-1}f(1)>0$.
    \end{mainresult}

    Bernstein's representation implies that every nonconstant classical completely monotone function on $(0,\infty)$ satisfies the strict difference condition when $D\ge2$. Thus Main Result~\ref{main:checkerboard} proves the checkerboard conjecture for this class on every product of even cycles. Strict decrease and strict convexity alone do not suffice, the example on $\Z_8^2$ in Section~\ref{sec:checkerboards} gives a different configuration of the same size with lower energy.

    The extension to arbitrary even side lengths addresses a specific obstruction to the earlier spectral argument. As observed in~\cite[Remark~5.5]{ref5}, when a side length is greater than two and congruent to two modulo four, the checkerboard frequency need not correspond to the smallest eigenvalue of the energy kernel. A proof must therefore use information beyond that spectral relaxation. We retain such information through the integrality of local occupancy counts. Its role is already visible for three points in $\Z_6$, the ordinary Fourier LP admits a feasible distribution with $b_1$-energy $3/2$, below the minimum energy $2$ of an actual set, whereas our integer box constraints exclude that distribution (Section~\ref{sec:checkerboards}).

    More precisely, we expand every binomial potential as a nonnegative combination of box autocorrelations. Each autocorrelation expresses a weighted pair count as a sum of squared occupancies over translations. For either checkerboard, each count is one of the two integers bracketing half the box volume, so all these sums attain their minimum permitted by integrality. This connection between energy and local counts has a counterpart in discrepancy identities for Hamming space~\cite{ref4}. The same method also gives universally optimal sets of densities $1/3$ and $2/3$ whenever every side length is divisible by three, including products with unequal side lengths (Theorem~\ref{cor:index-three}).

    Half the volume of an even torus is also the maximum size of a set with no pair at distance one. On $\Z_q^n$ this maximum is $M=\lfloor q/2\rfloor q^{n-1}$. At that cardinality the nearest neighbor potential is minimized by every maximum independent set. Our next question is which tori admit a maximum independent set that also minimizes every completely monotone potential.

    \begin{mainresult}[Theorem~\ref{thm:odd} and Corollary~\ref{cor:existence}]\label{main:classification}
        Let $q\ge2$, $n\ge1$, and $M=\lfloor q/2\rfloor q^{n-1}$. A universally optimal $M$-point set in $\Z_q^n$ exists exactly when $q$ is even, $q=3$, $n=1$, or $(q,n)=(5,2)$. The same parameter list characterizes the existence of a common minimizer for all exponential potentials $c^r$, $0<c<1$. In each of these cases the universally optimal sets are precisely the maximum independent sets. They form a single Lee isometry class and, for each fixed $c\in(0,1)$, are exactly the minimizers of the exponential energy.
    \end{mainresult}

    The obstruction in the remaining cases comes from the different behavior at the two ends of the exponential family. For sufficiently small $c$, every minimizing set must be independent. On the excluded odd tori every maximum independent set is isometric to a stripe by the classification of Arras and Joos~\cite{refAJ}. For $c$ sufficiently close to one, we construct a set with energy below that of the stripe. In the exceptional plane $\Z_5^2$, there are exactly ten optimal ten-point sets on the labeled torus, given by $x+y\in\{a,a+2\}$ or $x-y\in\{a,a+2\}$ for $a\in\Z_5$.

    For sets of density $1/5$ in $\Z_5^n$, we seek a criterion that characterizes all universally optimal configurations. Although the diameter grows with $n$, a test involving only the first two distances suffices. Let $\nu_r(C)$ be the number of unordered pairs in $C$ at Lee distance $r$.

    \begin{mainresult}[Theorem~\ref{thm:quinary-characterization}]\label{main:latin}
        For $n\ge2$, a $5^{n-1}$-point set $C\subseteq\Z_5^n$ is universally optimal if and only if it minimizes $(2n-1)\nu_1(C)+\nu_2(C)$ at that cardinality.
    \end{mainresult}

    The equality conditions in this two distance bound give an explicit restriction on the Fourier support of $1_C$. Besides characterizing optimal sets, this restriction yields constructions obtained by lifting Latin hypercubes (Theorem~\ref{thm:latin}).

    Section~\ref{sec:preliminaries} fixes notation and records the Fourier bound. Section~\ref{sec:integer-boxes} develops the integer box inequalities, which are applied to checkerboards and character fibers in Section~\ref{sec:checkerboards}. Section~\ref{sec:odd} establishes the classification at the independence number. Finally, Section~\ref{sec:quinary} proves the characterization at density $1/5$ in $\Z_5^n$ and gives the Latin lift construction.

\section{Preliminaries}
\label{sec:preliminaries}

Let
\[
G=\prod_{i=1}^n\Z_{q_i},\qquad q_i\ge2,\qquad
V=\prod_{i=1}^n q_i,\qquad
s_i=\lfloor q_i/2\rfloor,\qquad D=\sum_{i=1}^n s_i.
\]
The Lee metric depends on the specified cyclic decomposition. Write
\[
d(z)=\sum_{i=1}^n\ell_{q_i}(z_i),\qquad d_L(x,y)=d(x-y).
\]
Throughout, $n\ge1$ and $1\le M\le V$, and energy comparisons are
between sets of cardinality $M$. For $x\in C$, put
\[
U_f(x;C)=\sum_{y\in C\setminus\{x\}}f(d_L(x,y)),\qquad
E_f(C)=\frac1M\sum_{x\in C}U_f(x;C),\qquad
E_f^{\max}(C)=\max_{x\in C}U_f(x;C).
\]
The total energy of unordered pairs is $ME_f(C)/2$. If $C$ minimizes
$E_f$ and $U_f(x;C)$ is constant on $C$, then every $M$-point set $A$
satisfies
\[
E_f^{\max}(A)\ge E_f(A)\ge E_f(C)=E_f^{\max}(C).
\]
Thus a mean energy bound also gives a bound for the maximum individual
energy whenever the minimizing set has constant individual energy, as in~\cite[Section~2]{ref5}.

The distance distribution can be encoded by the normalized enumerator
\begin{equation}
W_C(z)=\frac1M\sum_{x,y\in C}z^{d_L(x,y)},\qquad
E_f(C)=\sum_{r=1}^D[z^r]W_C(z)f(r).
\label{eq:2.2}
\end{equation}
In particular, $[z^0]W_C(z)=1$ and $W_C(1)=M$. For a coset
$C=x_0+H$,
\[
W_C(z)=\sum_{h\in H}z^{d(h)}.
\]
We set $\binom ab=0$ for a nonnegative integer $a$ and an integer
$b\notin[0,a]$. A coefficient at a negative index is also understood
to be zero.

For potentials on a finite set of distances, we use the finite difference notion of complete monotonicity.
\begin{definition}\label{def:universal}
Write $\Delta f(t)=f(t+1)-f(t)$. A sequence $f$ on $\{1,\ldots,D\}$ is
\emph{completely monotone} if
\[
(-1)^k\Delta^k f(t)\ge0
\qquad(k\ge0,\ t\ge1,\ t+k\le D).
\]
A set is \emph{universally optimal} if it minimizes $E_f$ among sets of
the same cardinality for every such sequence.
\end{definition}

This definition follows~\cite{ref12}. A function
$f\in C^\infty(0,\infty)$ is \emph{classically completely monotone}
if $(-1)^k f^{(k)}(r)\ge0$ for all $k\ge0$ and $r>0$.
Successive applications of the fundamental theorem of calculus give
\begin{align*}
\Delta f(t)&=\int_0^1 f'(t+u)\,du,\\
(-1)^k\Delta^k f(t)
&=\int_{[0,1]^k}(-1)^k f^{(k)}(t+u_1+\cdots+u_k)
\,du_1\cdots du_k\ge0\qquad(k\ge1).
\end{align*}
Thus restricting a classically completely monotone function to the integer
distances gives a completely monotone sequence.

Cohn and Zhao's description of this cone~\cite[Lemma~4]{ref12}
provides a finite family of potentials that suffices to test universal
optimality. We state their representation on $\{1,\ldots,D\}$.
\begin{lemma}\label{lem:fundamental}
Put $b_j(t)=\binom{D-t}{j}$ for $0\le j\le D-1$. Every completely
monotone sequence has the unique expansion
\begin{equation}
f(t)=\sum_{j=0}^{D-1}c_jb_j(t),\qquad
c_j=(-1)^j\Delta^jf(D-j)\ge0.
\label{eq:2.3}
\end{equation}
Conversely, every nonnegative linear combination of the $b_j$ is
completely monotone.
\end{lemma}

For Fourier calculations, we extend each fundamental potential to the
diagonal by $b_j(0)=\binom Dj$, these values are excluded from the
pair energy. Linearity gives
\[
E_f(A)-E_f(C)
=\sum_{j=0}^{D-1}c_j\bigl(E_{b_j}(A)-E_{b_j}(C)\bigr).
\]
Consequently, a set is universally optimal if and only if it minimizes
all $D$ fundamental energies. The case $j=0$ is automatic, since
$E_{b_0}(C)=M-1$.

We next formulate a Fourier linear program for these energy minimization
problems. Identify the dual group with $G$ by
\[
\chi_u(z)=\exp\left(2\pi i\sum_{i=1}^n\frac{u_i z_i}{q_i}\right),
\qquad
\widehat h(u)=\frac1V\sum_{z\in G}h(z)\overline{\chi_u(z)}.
\]
Character orthogonality and Fourier inversion take the forms
\[
\sum_{u\in G}\chi_u(z)=V\mathbf1_{\{z=0\}},\qquad
h(z)=\sum_{u\in G}\widehat h(u)\chi_u(z).
\]
All auxiliary kernels below are real and even. For a set $C$, define its
character sums and normalized difference distribution by
\[
S_C(u)=\sum_{x\in C}\chi_u(x),\qquad
a_C(z)=\frac{|C\cap(C+z)|}{M}.
\]
The quantity $Ma_C(z)$ counts the ordered pairs $(x,y)\in C^2$
with $x-y=z$. Therefore
\[
a_C(0)=1,\qquad a_C(z)=a_C(-z),\qquad
\sum_{z\in G}a_C(z)=M.
\]
Character orthogonality gives
\begin{align*}
Ta_C(u)&:=\sum_{z\in G}a_C(z)\chi_u(z)
=\frac1M\sum_{x,y\in C}\chi_u(x-y)
=\frac{|S_C(u)|^2}{M},\\
\sum_{u\in G}|S_C(u)|^2
&=\sum_{x,y\in C}\sum_{u\in G}\chi_u(x-y)
=\sum_{x,y\in C}V\mathbf1_{\{x=y\}}=VM.
\end{align*}
Fourier inversion expresses the average kernel value in terms of the
character sums
\[
\sum_{z\in G}a_C(z)h(z)
=\frac1M\sum_{x,y\in C}h(x-y)
=\frac1M\sum_{u\in G}\widehat h(u)|S_C(u)|^2.
\]
We use this normalization throughout the Fourier bounds below.

Relaxing the requirement that $a$ arise from an actual set leads to the
ordinary Fourier linear program (LP). At cardinality $M$, it minimizes
$\sum_{z\ne0}a(z)f(d(z))$ subject to
\begin{equation}
a(0)=1,\quad a(z)=a(-z)\ge0,\quad
\sum_{z\in G}a(z)=M,\quad
Ta(u):=\sum_{z\in G}a(z)\chi_u(z)\ge0\quad(u\in G).
\label{eq:2.5}
\end{equation}
These are Delsarte's inner distribution inequalities for the translation
scheme of $G$~\cite[Theorem~3.3]{ref14}. We call $C$
\emph{ordinary Fourier LP universally optimal} if its difference
distribution $a_C$ attains this LP minimum for every completely
monotone potential. This property implies universal optimality because
every actual set has a feasible difference distribution.

We record the corresponding dual bound, together with the equality
conditions that will be used later.
\begin{theorem}\label{thm:fourier}
Let $h:G\to\R$ be an even function satisfying
\[
h(z)\le f(d(z))\quad(z\ne0),\qquad
\widehat h(u)\ge0\quad(u\ne0).
\]
Every vector feasible for \eqref{eq:2.5}, and hence every $M$-point
set, satisfies
\[
E_f\ge M\widehat h(0)-h(0).
\]
An $M$-point set $C$ attains the bound if and only if
\[
h(x-y)=f(d_L(x,y))\quad(x\ne y\in C),\qquad
\widehat h(u)|S_C(u)|^2=0\quad(u\ne0).
\]
\end{theorem}

\begin{proof}
For a feasible vector $a$, set
$E_f(a)=\sum_{z\ne0}a(z)f(d(z))$. Since $a(0)=1$ and $Ta(0)=M$,
\begin{align*}
\sum_{z\in G}a(z)h(z)
&=\sum_{u\in G}\widehat h(u)Ta(u)
=M\widehat h(0)+\sum_{u\ne0}\widehat h(u)Ta(u),\\
E_f(a)-\bigl(M\widehat h(0)-h(0)\bigr)
&=\sum_{z\ne0}a(z)\bigl(f(d(z))-h(z)\bigr)
  +\sum_{u\ne0}\widehat h(u)Ta(u)\ge0.
\end{align*}
For an actual set $C$, substituting $a=a_C$ gives the exact gap identity
\begin{equation}
\begin{aligned}
E_f(C)-\bigl(M\widehat h(0)-h(0)\bigr)
=\frac1M\sum_{\substack{x,y\in C\\x\ne y}}
\bigl(f(d_L(x,y))-h(x-y)\bigr)
+\frac1M\sum_{u\ne0}\widehat h(u)|S_C(u)|^2.
\end{aligned}
\label{eq:fourier-gap}
\end{equation}
Since every term is nonnegative, equality holds precisely when every
summand in both sums vanishes. These are exactly the stated conditions.
\end{proof}

To compare configurations up to equivalence, we also record the form
of the ambient Lee isometries. For $q\ge3$, the Lee distance on $\Z_q^n$ is the graph distance
on the Cartesian power of $C_q$. Hence its isometries are precisely
the graph automorphisms. For $q\ge5$, the cycle $C_q$ is prime for
the Cartesian product, a product of two nontrivial connected graphs
contains a $4$-cycle, whereas $C_q$ does not. The product automorphism
theorem~\cite[Theorem~6.10]{HIK}, together with
$\operatorname{Aut}(C_q)=\{t\mapsto a\pm t:a\in\Z_q\}$, gives
the following specialization.

\begin{lemma}\label{lem:isometries}
For $q\ge5$, every Lee isometry of $\Z_q^n$ has the form
\[
(x_1,\ldots,x_n)\longmapsto
(a_1+\varepsilon_1x_{\sigma(1)},\ldots,
 a_n+\varepsilon_nx_{\sigma(n)}),
\]
where $a\in\Z_q^n$, $\sigma$ is a permutation of the coordinates,
and $\varepsilon_i\in\{\pm1\}$. Every such map sends subgroup
cosets to subgroup cosets.
\end{lemma}

    \section{Box-strengthened LP bound}
    \label{sec:integer-boxes}

    We first express each binomial kernel as a nonnegative combination of box autocorrelations. Each autocorrelation turns a weighted pair count into a sum of squared occupancies over all translates. The integrality of these occupancies then yields a refinement of the ordinary Fourier bound.

    For $s=\lfloor q/2\rfloor$ and $1\le a\le s$, let
    $I_a=\{0,\ldots,a-1\}\subseteq\Z_q$. Its autocorrelation is
    \[
    K_a(z)=|I_a\cap(I_a+z)|=(a-\ell_q(z))_+,
    \qquad x_+=\max(x,0).
    \]
    Indeed, an integer difference $h\in\{1-a,\ldots,a-1\}$ occurs
    $a-|h|$ times among pairs in $I_a^2$. Since $2a-2<q$, these
    differences are distinct modulo $q$, and hence
    \[
    K_a(z)=\sum_{h=1-a}^{a-1}(a-|h|)\mathbf1_{\{h\equiv z\pmod q\}}
    =\begin{cases}
      a-\ell_q(z),&\ell_q(z)<a,\\
      0,&\ell_q(z)\ge a.
    \end{cases}
    \]
    We also put $I_q=\Z_q$, so $K_q(z)=q$ for every $z$.
    An \emph{allowed box} is a product of cyclic intervals whose
    lengths in coordinate $i$ belong to $\{1,\ldots,s_i,q_i\}$.
    The next lemma gives the required nonnegative expansion, first on a
    cyclic factor and then on the product group.

    \begin{lemma}
        \label{lem:box-expansion}
        On $\Z_q$, with $s=\lfloor q/2\rfloor$, the generating function of the binomial kernels satisfies
        \begin{equation}
            (1+t)^{s-\ell_q(z)}
            =\frac{K_q(z)}q+tK_s(z)
            +t^2\sum_{a=1}^{s-1}(1+t)^{s-a-1}K_a(z).
            \label{eq:3.3}
        \end{equation}
        For a length vector $\boldsymbol{a}\in\mathcal A:=\prod_i\{1,\ldots,s_i,q_i\}$, define
        \[
        I_{\boldsymbol{a}}=\prod_i I_{a_i},\qquad
        Q_{\boldsymbol{a}}=|I_{\boldsymbol{a}}|,\qquad
        K_{\boldsymbol{a}}(z)=|I_{\boldsymbol{a}}\cap(I_{\boldsymbol{a}}+z)|.
        \]
        Set
        \[
        \alpha_{i,q_i}(t)=q_i^{-1},\qquad
        \alpha_{i,s_i}(t)=t,\qquad
        \alpha_{i,a}(t)=t^2(1+t)^{s_i-a-1}\quad(1\le a<s_i),
        \]
        \[
        w_{j,\boldsymbol{a}}=[t^j]\prod_i\alpha_{i,a_i}(t).
        \]
        Then $w_{j,\boldsymbol{a}}\ge0$ and
        \begin{equation}
            \binom{D-d(z)}j
            =\sum_{\boldsymbol{a}\in\mathcal A}w_{j,\boldsymbol{a}}K_{\boldsymbol{a}}(z)
            \qquad (0\le j\le D).
            \label{eq:3.4}
        \end{equation}
        For $j\le D-1$, the left-hand side is $b_j(d(z))$, and for $j=D$, it is the indicator of the origin.
    \end{lemma}

    \begin{proof}
    Put $r=\ell_q(z)$ and regard $t$ as a formal variable. The
    generating function for subsets of $S_r=\{r+1,\ldots,s\}$ is
    \[
    \sum_{A\subseteq S_r}t^{|A|}=(1+t)^{s-r}.
    \]
    The subsets of sizes zero and one contribute $1+(s-r)t$.
    Every subset of size at least two has the unique form
    \[
    A=\{m,a+1\}\cup B,\qquad
    r+1\le m\le a\le s-1,\qquad
    B\subseteq\{a+2,\ldots,s\},
    \]
    where $a+1$ is its second smallest element. Consequently,
    \begin{align*}
    (1+t)^{s-r}
    &=1+(s-r)t+
      \sum_{a=r+1}^{s-1}\sum_{m=r+1}^{a}
      \sum_{B\subseteq\{a+2,\ldots,s\}}t^{2+|B|}\\
    &=1+(s-r)t+
      t^2\sum_{a=r+1}^{s-1}(a-r)(1+t)^{s-a-1}\\
    &=\frac{K_q(z)}q+tK_s(z)+
      t^2\sum_{a=1}^{s-1}(1+t)^{s-a-1}K_a(z).
    \end{align*}
    In the last line we used $K_s(z)=s-r$ and
    $K_a(z)=(a-r)_+$. Empty sums are zero, so the calculation also
    covers $r=s$ and $s=1$. At fixed subset size, this partition
    gives the binomial convolution in~\cite[Eq.~(5.26)]{GKP}.

    For the product group, denote the autocorrelation in coordinate
    $i$ by $K^{(i)}_{a_i}$. Cartesian products give
    \[
    K_{\boldsymbol a}(z)=\prod_{i=1}^nK^{(i)}_{a_i}(z_i),
    \qquad
    D-d(z)=\sum_{i=1}^n\bigl(s_i-\ell_{q_i}(z_i)\bigr).
    \]
    Applying the one-dimensional identity in each coordinate yields
    \begin{align*}
    (1+t)^{D-d(z)}
    &=\prod_{i=1}^n
      \left(\sum_{a_i\in\{1,\ldots,s_i,q_i\}}
      \alpha_{i,a_i}(t)K^{(i)}_{a_i}(z_i)\right)\\
    &=\sum_{\boldsymbol a\in\mathcal A}
      \left(\prod_{i=1}^n\alpha_{i,a_i}(t)\right)
      K_{\boldsymbol a}(z).
    \end{align*}
    Coefficient extraction now gives
    \[
    \binom{D-d(z)}j
    =[t^j](1+t)^{D-d(z)}
    =\sum_{\boldsymbol a\in\mathcal A}
      w_{j,\boldsymbol a}K_{\boldsymbol a}(z).
    \]
    To verify the nonnegativity of the coefficients explicitly, set
    \[
    F=\{i:a_i=q_i\},\qquad
    T=\{i:a_i=s_i\},\qquad
    S=\{i:a_i<s_i\}.
    \]
    These three sets partition the coordinates. With
    $r_{\boldsymbol a}=|T|+2|S|$ and
    $h_{\boldsymbol a}=\sum_{i\in S}(s_i-a_i-1)$, we have
    \[
    \prod_i\alpha_{i,a_i}(t)
    =\left(\prod_{i\in F}q_i^{-1}\right)
       t^{r_{\boldsymbol a}}(1+t)^{h_{\boldsymbol a}},
    \qquad
    w_{j,\boldsymbol a}
    =\left(\prod_{i\in F}q_i^{-1}\right)
       \binom{h_{\boldsymbol a}}{j-r_{\boldsymbol a}}\ge0.
    \]
    Finally,
    \[
    \binom{D-d(z)}D=\mathbf1_{\{d(z)=0\}}
    =\mathbf1_{\{z=0\}},
    \]
    which proves the assertion for $j=D$.
\end{proof}

    For any box $B\subseteq G$, its autocorrelation $K_B(z)=|B\cap(B+z)|$ satisfies
    \[
    \widehat K_B(u)
    =\frac1V\sum_{x,y\in B}\overline{\chi_u(x-y)}
    =\frac1V\left|\sum_{x\in B}\chi_u(x)\right|^2\ge0.
    \]
    Together with Lemma~\ref{lem:box-expansion}, this shows that all fundamental kernels have nonnegative Fourier coefficients, regardless of the side lengths. To obtain a stronger bound, we now use the fact that the occupancy of each translated box is an integer.

    Fix a box $B$ of volume $Q$ and an $M$-point set $C$. For $t\in G$
    put $N_t=|C\cap(B+t)|=\sum_{x\in C}\mathbf1_B(x-t)$.
    Counting incidences first by points and then by ordered pairs gives
    \begin{equation}
    	\begin{split}
    		&\sum_{t\in G}N_t
    		=\sum_{x\in C}\sum_{t\in G}\mathbf1_B(x-t)=MQ,\\
    		&\sum_{t\in G}N_t^2
    		=\sum_{x,y\in C}\sum_{t\in G}
    		\mathbf1_B(x-t)\mathbf1_B(y-t)
    		=\sum_{x,y\in C}K_B(x-y)
    		=M\sum_{z\in G}a_C(z)K_B(z).
    	\end{split}
    	\label{eq:3.5}
    \end{equation}
    Write
    \[
    \lambda=\frac MV,\qquad
    \lambda Q=k+\theta,\qquad
    k=\lfloor\lambda Q\rfloor,\qquad 0\le\theta<1.
    \]
    Since $(r-k)(r-k-1)\ge0$ for every integer $r$,
    \begin{align*}
    0
    &\le\sum_t(N_t-k)(N_t-k-1)\\
    &=\sum_tN_t^2-(2k+1)\sum_tN_t+Vk(k+1)\\
    &=M\sum_z a_C(z)K_B(z)-(2k+1)MQ+Vk(k+1).
    \end{align*}
    Dividing by $M$ proves
    \begin{equation}
    \sum_z a_C(z)K_B(z)\ge(2k+1)Q-\frac VM k(k+1),
    \qquad k=\left\lfloor\frac{MQ}V\right\rfloor.
    \label{eq:3.7}
    \end{equation}
    Equality holds if and only if $N_t\in\{k,k+1\}$ for every $t$.
    If $\theta=0$, then
    \[
    \sum_t(N_t-k)=MQ-Vk=0,\qquad N_t-k\in\{0,1\},
    \]
    so every $N_t=k$. Thus in all cases equality is equivalent to
    \[
    N_t\in\{\lfloor\lambda Q\rfloor,\lceil\lambda Q\rceil\}
    \qquad(t\in G).
    \]
    This is the minimum sum of squares property of balanced integer
    lists~\cite[Lemmas~2.8 and~2.15]{BCL}.

    Since $K_B(0)=Q$, the corresponding off-diagonal energy satisfies
    \begin{align*}
    E_{K_B}(C)
    &:=\sum_{z\ne0}a_C(z)K_B(z)\\
      &\ge 2kQ-\frac1\lambda k(k+1)\\
    &=\frac{2k(k+\theta)-k(k+1)}\lambda\\
    &=\lambda Q^2-Q+\frac{\theta(1-\theta)}\lambda.
    \end{align*}
    The kernel $K_B$ need not depend only on Lee distance, but the
    ordinary Fourier argument still gives the bound
    \[
    M\widehat K_B(0)-K_B(0)
    =\frac{MQ^2}{V}-Q=\lambda Q^2-Q.
    \]
    The integer correction $\theta(1-\theta)/\lambda$ is positive
    precisely when $MQ/V\notin\Z$.

    The \emph{box-strengthened LP} is \eqref{eq:2.5} with \eqref{eq:3.7} imposed for every
    $B=I_{\boldsymbol{a}}$, $\boldsymbol{a}\in\mathcal A$. This is an explicit list of
    $\prod_i(s_i+1)$ inequalities. Translating a box does not change its autocorrelation, so no
    further LP constraints are needed for its translates.
    Since the expansion in Lemma~\ref{lem:box-expansion} has nonnegative
    coefficients, simultaneous equality in the box inequalities gives
    equality for every completely monotone potential. The following
    theorem makes this criterion precise.

    \begin{theorem}
        \label{thm:box}
        Suppose that an $M$-point set $C\subseteq G$ satisfies
        \begin{equation}
            |C\cap(B+t)|\in
            \left\{\left\lfloor\frac{M|B|}V\right\rfloor,
            \left\lceil\frac{M|B|}V\right\rceil\right\}
            \label{eq:3.9}
        \end{equation}
        for every allowed box $B$ and every $t\in G$. Then $C$ attains the box-strengthened LP for every completely monotone potential.

        More precisely, write $f=\sum_j c_j b_j$, let $\eta_{f,\boldsymbol a}=\sum_j c_j w_{j,\boldsymbol a}$, and put $k_{\boldsymbol a}=\lfloor M Q_{\boldsymbol a}/V\rfloor$. Every vector feasible for the box-strengthened LP satisfies
        \begin{equation}
            E_f\ge
            \sum_{\boldsymbol{a}\in\mathcal A}\eta_{f,\boldsymbol{a}}
            \left(2k_{\boldsymbol{a}}Q_{\boldsymbol{a}}
            -\frac VM k_{\boldsymbol{a}}(k_{\boldsymbol{a}}+1)\right).
            \label{eq:3.11}
        \end{equation}
        An actual set attains this bound if and only if \eqref{eq:3.9} holds for each box with
        $\eta_{f,\boldsymbol{a}}>0$.
    \end{theorem}

    \begin{proof}
    By Lemmas~\ref{lem:fundamental} and~\ref{lem:box-expansion},
    \begin{align*}
    f(d(z))
    &=\sum_{j=0}^{D-1}c_j\sum_{\boldsymbol a\in\mathcal A}
      w_{j,\boldsymbol a}K_{\boldsymbol a}(z)\\
    &=\sum_{\boldsymbol a\in\mathcal A}
      \eta_{f,\boldsymbol a}K_{\boldsymbol a}(z)
      \qquad(z\ne0),
    \end{align*}
    where $\eta_{f,\boldsymbol a}=\sum_jc_jw_{j,\boldsymbol a}\ge0$.
    Define
    \[
    L_{\boldsymbol a}
    =2k_{\boldsymbol a}Q_{\boldsymbol a}
     -\frac VM k_{\boldsymbol a}(k_{\boldsymbol a}+1),
    \qquad
    L_f=\sum_{\boldsymbol a}\eta_{f,\boldsymbol a}L_{\boldsymbol a}.
    \]
    For every feasible vector $a$ of the strengthened LP,
    \begin{align*}
    E_f(a)
    &=\sum_{z\ne0}a(z)f(d(z))\\
    &=\sum_{\boldsymbol a}\eta_{f,\boldsymbol a}
      \left(\sum_z a(z)K_{\boldsymbol a}(z)-Q_{\boldsymbol a}\right)\\
    &\ge\sum_{\boldsymbol a}\eta_{f,\boldsymbol a}L_{\boldsymbol a}
      =L_f.
    \end{align*}
    The second equality uses $a(0)=1$ and
    $K_{\boldsymbol a}(0)=Q_{\boldsymbol a}$; the inequality is
    the nonnegative weighted sum of \eqref{eq:3.7}.

    To determine when equality holds, let $A$ be an actual $M$-point set and put
    $N_{\boldsymbol a,t}=|A\cap(I_{\boldsymbol a}+t)|$.
    By \eqref{eq:3.5}, the slack for one box shape is
    \begin{align*}
    &\sum_z a_A(z)K_{\boldsymbol a}(z)
      -Q_{\boldsymbol a}-L_{\boldsymbol a}\\
    =&\frac1M\left(
      \sum_tN_{\boldsymbol a,t}^2
      -(2k_{\boldsymbol a}+1)MQ_{\boldsymbol a}
      +Vk_{\boldsymbol a}(k_{\boldsymbol a}+1)\right)\\
    =&\frac1M\sum_t
      (N_{\boldsymbol a,t}-k_{\boldsymbol a})
      (N_{\boldsymbol a,t}-k_{\boldsymbol a}-1).
    \end{align*}
    Summing over the box shapes gives the exact energy gap
    \begin{equation}
    \begin{split}
    E_f(A)-L_f
    =\frac1M\sum_{\boldsymbol a}\eta_{f,\boldsymbol a}
       \sum_t(N_{\boldsymbol a,t}-k_{\boldsymbol a})
       (N_{\boldsymbol a,t}-k_{\boldsymbol a}-1).
    \end{split}
    \label{eq:box-energy-gap}
    \end{equation}
    Each product in the inner sum is nonnegative. Equality
    holds precisely when, for every $t\in G$ and every shape with
    $\eta_{f,\boldsymbol a}>0$,
    \begin{align*}
    (N_{\boldsymbol a,t}-k_{\boldsymbol a})
       (N_{\boldsymbol a,t}-k_{\boldsymbol a}-1)=0\qquad\Longleftrightarrow\quad
       N_{\boldsymbol a,t}\in\{k_{\boldsymbol a},k_{\boldsymbol a}+1\}.
    \end{align*}
    If $MQ_{\boldsymbol a}/V$ is nonintegral, these are exactly
    its floor and ceiling. In the integral mean case, if
    $MQ_{\boldsymbol a}/V=k_{\boldsymbol a}$, then
    \[
    \sum_t(N_{\boldsymbol a,t}-k_{\boldsymbol a})=0,
    \qquad N_{\boldsymbol a,t}-k_{\boldsymbol a}\ge0,
    \]
    forcing every summand to vanish. Thus equality is equivalent to
    \eqref{eq:3.9} for every box shape with positive weight.

    Finally, every actual set has a difference distribution feasible
    for the strengthened LP. The balance condition for $C$ makes
    every term in \eqref{eq:box-energy-gap} zero, so
    \[
    E_f(C)=L_f
    =\min_{a\text{ feasible}}E_f(a)
    \le E_f(A)\qquad(|A|=M).
    \]
    This proves both LP attainment and universal optimality.
\end{proof}

    The box criterion is preserved under taking complements. Indeed, for $Q=|B|$,
    \[
    |(G\setminus C)\cap(B+t)|=Q-|C\cap(B+t)|,
    \]
    and the rounding identities
    \[
    Q-\lfloor\lambda Q\rfloor=\lceil(1-\lambda)Q\rceil,\qquad
    Q-\lceil\lambda Q\rceil=\lfloor(1-\lambda)Q\rfloor
    \]
    give \eqref{eq:3.9} at the complementary density.

    \section{Universally optimal sets with size $\frac{|G|}{2}$ and $\frac{|G|}{3}$}
    \label{sec:checkerboards}

    We apply the box criterion to fibers of characters of orders two and three. For the parity classes, we also obtain uniqueness and optimality for maximum individual energy. A union of two fibers in $\Z_5^2$ will supply another positive example needed for the classification in Section~\ref{sec:odd}.

    \begin{theorem}
        \label{thm:checkerboard}
        Suppose that every $q_i$ is even. The two parity classes
        \[
        H_0=\left\{x\in G:\sum_i x_i\equiv0\pmod2\right\},
        \qquad H_1:=G\setminus H_0
        \]
        both have size $V/2$. For every completely monotone $f$, they attain the box-strengthened LP and minimize the mean, total pair, and maximum individual energies.
        If $D\ge2$ and
        \begin{equation}
            (-1)^{D-1}\Delta^{D-1}f(1)>0,
            \label{eq:4.1}
        \end{equation}
        they are the only $V/2$-point minimizers of the mean, total pair, and maximum individual energies. If $D=1$, the two parity classes
        are the only competitors.
    \end{theorem}

    \begin{proof}
    Because every $q_i$ is even, the map
    $\pi(x)=\sum_i x_i\pmod2$ is a well-defined surjective
    homomorphism. Write $H=\ker\pi$ and $H_\varepsilon=\pi^{-1}(\varepsilon)$. Then
    \[
    |H_0|=|H_1|=\frac V2,\qquad
    \mathbf1_{H_\varepsilon}(x)
    =\frac{1+(-1)^\varepsilon\prod_i(-1)^{x_i}}2.
    \]
    For a translated box $B=\prod_iJ_i$, write
    $J_i=t_i+\{0,\ldots,a_i-1\}$ and $Q=\prod_i a_i$.
    The box occupancies are determined by the alternating sums
    \[
    S_i:=\sum_{u\in J_i}(-1)^u
    =(-1)^{t_i}\sum_{j=0}^{a_i-1}(-1)^j
    =(-1)^{t_i}\frac{1-(-1)^{a_i}}2.
    \]
    Hence $|H_\varepsilon\cap B|
    =\frac12\left(Q+(-1)^\varepsilon\prod_iS_i\right)$, where
    \begin{align*}
    \prod_iS_i&=
    \begin{cases}
    0,&Q\text{ even},\\
    (-1)^{\sum_i t_i},&Q\text{ odd},
    \end{cases}
    \end{align*}
    and therefore $|H_\varepsilon\cap B|
    \in\{\lfloor Q/2\rfloor,\lceil Q/2\rceil\}.$
    Theorem~\ref{thm:box} proves LP attainment and minimization of
    the mean and total pair energies.

    The subgroup structure also makes the individual energies constant on each parity class.
    For $x\in H_\varepsilon$, the map $y\mapsto y-x$ identifies
    $H_\varepsilon\setminus\{x\}$ with $H\setminus\{0\}$, so
    \[
    U_f(x;H_\varepsilon)=\sum_{h\in H\setminus\{0\}}f(d(h))
    =E_f(H)=E_f^{\max}(H_\varepsilon).
    \]
    Consequently, every $V/2$-point set $A$ satisfies $E_f^{\max}(A)\ge E_f(A)\ge E_f(H)
    =E_f^{\max}(H_\varepsilon).$

    For uniqueness, assume $D\ge2$ and \eqref{eq:4.1}. In the expansion
    \eqref{eq:2.3},
    \[
    c_{D-1}=(-1)^{D-1}\Delta^{D-1}f(1)>0,
    \qquad
    b_{D-1}(r)=\mathbf1_{\{r=1\}}\quad(1\le r\le D).
    \]
    Let $e(A)$ count the unordered distance one pairs in $A$.
    Since $e(H)=0$ and $H$ minimizes every fundamental energy,
    \begin{align*}
    E_f(A)-E_f(H)
    =\sum_{j=0}^{D-1}c_j\bigl(E_{b_j}(A)-E_{b_j}(H)\bigr)
    \ge c_{D-1}\bigl(E_{b_{D-1}}(A)-E_{b_{D-1}}(H)\bigr)
      =\frac{2c_{D-1}}{|A|}e(A).
    \end{align*}
    Thus any mean energy minimizer is independent in the
    distance one graph $\Gamma$. This graph is regular with degree $R=\sum_{i=1}^n\bigl(2-\mathbf1_{\{q_i=2\}}\bigr),$ and its total number of edges is $|E(\Gamma)|=\frac{RV}{2}$.
    For an independent $A$ of size $V/2$, the number of edges
    crossing from $A$ to its complement is
    \[
    e(A,G\setminus A)=R|A|=\frac{RV}{2}=|E(\Gamma)|.
    \]
    Hence the complement of $A$ is independent as well. In terms of
    $g(x)=2\mathbf1_A(x)-1$, this says $g(x+e_i)=-g(x)\ (x\in G,\ 1\le i\le n).$
    Iterating along coordinate steps yields
    \[
    g(x)=g(0)(-1)^{\sum_i x_i},\qquad g(0)\in\{-1,1\},
    \]
    so $A=H_0$ or $A=H_1$. If $A$ minimizes the maximum
    individual energy, the chain
    \[
    E_f(H)\le E_f(A)\le E_f^{\max}(A)=E_f^{\max}(H)=E_f(H)
    \]
    gives the same conclusion. Finally,
    $D=1$ and even side lengths imply $n=1$, $q_1=2$.
    The two singleton sets are then the only competitors.
\end{proof}

    The strict difference condition in Theorem~\ref{thm:checkerboard} holds
    for every nonconstant classical completely monotone function on $(0,\infty)$.
    Indeed, Bernstein's theorem gives
    \[
    f(r)=\int_{[0,\infty)}e^{-rt}\,d\mu(t)
    \]
    for a positive measure $\mu$, as in the proof of
    \cite[Proposition~5.4]{ref5}. For $D\ge2$,
    \begin{align*}
    (-1)^{D-1}\Delta^{D-1}f(1)
    &=\sum_{j=0}^{D-1}(-1)^j\binom{D-1}{j}f(1+j)\\
    &=\int_{[0,\infty)}e^{-t}
      \sum_{j=0}^{D-1}\binom{D-1}{j}(-e^{-t})^j\,d\mu(t)\\
    &=\int_{[0,\infty)}e^{-t}(1-e^{-t})^{D-1}\,d\mu(t).
    \end{align*}
    The integrand is positive for $t>0$ and bounded above by
    $e^{-t}$, whose integral is $f(1)<\infty$. A nonconstant $f$
    has $\mu((0,\infty))>0$, so the last integral is positive.
    Theorem~\ref{thm:checkerboard} therefore gives uniqueness for
    every such $f$. For a constant $f\equiv a$, every $M$-point set
    has energy $E_f=(M-1)a$.

    Strict decrease and strict discrete convexity alone do not suffice
    for checkerboard optimality. On $\Z_8^2$, put
    \[
    f_0(r)=(3-r)_+,\quad
    C=\{x:x_1+x_2\equiv0\text{ or }1\pmod4\},\quad
    H=\{x:x_1+x_2\equiv0\pmod2\}.
    \]
    The four fibers of $x_1+x_2\pmod4$ have size $16$, so
    $|C|=|H|=32$. At a point of $C$ with coordinate sum $0\pmod4$,
    the unit steps remaining in $C$ are $e_1,e_2$. Among the eight
    distance two displacements $\{\pm2e_1,\pm2e_2,\pm e_1\pm e_2\},$
    exactly $\pm(e_1-e_2)$ remain in $C$. The same counts hold in the other
    occupied residue class, with the unit-step directions reversed.
    All eight distance two displacements remain in $H$. Thus, with $N_r(x;A)=|\{y\in A:d_L(x,y)=r\}|$,
    \[
    (N_1(x;C),N_2(x;C))=(2,2),\qquad
    (N_1(y;H),N_2(y;H))=(0,8),
    \]
    for all $x\in C$ and $y\in H$, and hence
    \[
    U_{f_0}(x;C)=2\cdot2+2\cdot1=6,
    \qquad U_{f_0}(y;H)=8.
    \]
    Now take $f(r)=f_0(r)+2^{-r}/16$. Its differences satisfy
    \begin{align*}
    \Delta f(r)&=-\mathbf1_{\{r\in\{1,2\}\}}
       -\frac{2^{-r-1}}{16}<0,\\
    \Delta^2f(r)&=\mathbf1_{\{r=2\}}
       +\frac{2^{-r-2}}{16}>0
    \end{align*}
    whenever the indicated differences are defined. Nevertheless,
    for every $x\in C$ and $y\in H$,
    \[
    U_f(x;C)\le6+31\frac{2^{-1}}{16}
      =6+\frac{31}{32}<8\le U_f(y;H).
    \]
    Thus $C$ has both smaller mean energy and smaller maximum
    individual energy than a checkerboard.

    The integer box constraints can strengthen the ordinary Fourier LP
    even at the checkerboard cardinality. This already occurs for three
    points in $\Z_6$. For $b_1(r)=3-r$, the checkerboard
    $H=\{0,2,4\}$ has $E_{b_1}(H)=2$. Consider the vector
    \[
    a=(1,\frac{1}{4},\frac{1}{4},1,\frac{1}{4},\frac{1}{4})
     =\tfrac14\mathbf1+\tfrac34
       (\mathbf1_{\{0\}}+\mathbf1_{\{3\}}).
    \]
    It satisfies
    \[
    a(0)=1,\quad a(z)=a(-z)\ge0,\quad
    Ta(0)=\sum_z a(z)=3,\quad
    Ta(u)=\tfrac34(1+(-1)^u)\ge0\quad(u\ne0).
    \]
    Hence $a$ is feasible for \eqref{eq:2.5}, but
    \[
    E_{b_1}(a)=\sum_{z\ne0}a(z)b_1(d(z))=2\left(\tfrac14\cdot2+\tfrac14\cdot1\right)
       =\frac32<2.
    \]
    For $I=\{0,1,2\}$, however,
    \[
    K_I=(3,2,1,0,1,2),\qquad
    \sum_z a(z)K_I(z)=3+2\left(\tfrac14\cdot2+\tfrac14\cdot1\right)
       =\frac92.
    \]
    Here $M=Q=3$, $V=6$, and $k=\lfloor MQ/V\rfloor=1$.
    The integer box bound requires
    \[
    \sum_z a(z)K_I(z)\ge(2k+1)Q-\frac VMk(k+1)
       =9-4=5>\frac92,
    \]
    so this vector is excluded by the strengthened LP.

    A similar box count argument applies to characters of order three,
    with interval sums of third roots of unity replacing the alternating sums.

    \begin{theorem}
        \label{cor:index-three}
        If every $q_i$ is divisible by three, then
        \[
        C_3(G)=\left\{x\in G:\sum_i x_i\equiv0\pmod3\right\}
        \]
        and its complement attain the box-strengthened LP for every completely monotone potential, at sizes $V/3$ and $2V/3$, respectively.
    \end{theorem}

    \begin{proof}
    The homomorphism $x\mapsto\sum_i x_i\pmod3$ is well defined
    and surjective, so its kernel $C_3(G)$ has size $V/3$.
    Set $\omega=e^{2\pi i/3}$. For a cyclic interval
    $J=b+\{0,\ldots,L-1\}$, the relation $\omega^{q_i}=1$
    permits the calculation
    \[
    U(J):=\sum_{x\in J}\omega^x
    =\omega^b\sum_{j=0}^{L-1}\omega^j
    =\begin{cases}
      0,&L\equiv0\pmod3,\\
      \omega^b,&L\equiv1\pmod3,\\
      -\omega^{b+2},&L\equiv2\pmod3.
    \end{cases}
    \]
    In particular, $|U(J)|\le1$. For $B=\prod_iJ_i$ and
    $Q=|B|$, the character identity
    $\mathbf1_{\{v\equiv0\pmod3\}}=(1+\omega^v+\omega^{-v})/3$
    gives
    \begin{align*}
    |C_3(G)\cap B|
    &=\frac13\sum_{x\in B}
      \left(1+\omega^{\sum_i x_i}+\omega^{-\sum_i x_i}\right)\\
    &=\frac13\left(Q+\prod_iU(J_i)+\overline{\prod_iU(J_i)}\right),\\
    \left|\,|C_3(G)\cap B|-\frac Q3\right|
    &\le\frac23\prod_i|U(J_i)|\le\frac23.
    \end{align*}
    If $Q=3h+v$, $v\in\{0,1,2\}$, the integer occupancy $N$
    therefore satisfies
    \[
    |3N-Q|\le2
    \quad\Longrightarrow\quad
    N\in\begin{cases}
       \{h\},&v=0,\\
       \{h,h+1\},&v=1\text{ or }2.
       \end{cases}.
    \]
    The complementary occupancy is $Q-N$ and hence lies in
    $\{\lfloor2Q/3\rfloor,\lceil2Q/3\rceil\}$.
    Theorem~\ref{thm:box} applies at both cardinalities.
\end{proof}

    When every $q_i=3$, the coordinate sum code is the ternary single parity check
    code in Hamming space. Its universal optimality is also the MDS case
    of~\cite[Table~I]{ref12} for $n\ge2$ and is immediate for $n=1$.
    This supplies the ternary case of Main Result~\ref{main:classification}.
    To treat the exceptional parameter pair $(q,n)=(5,2)$, we use a union
    of two fibers of the coordinate sum.

    \begin{proposition}
        \label{prop:ten}
        The set $C_{10}=\{(x,y)\in\Z_5^2:x+y\in\{0,2\}\}$
        and its complement attain the box-strengthened LP for every completely monotone potential at sizes 10 and 15, respectively.
    \end{proposition}

    \begin{proof}
    Each fiber of $(x,y)\mapsto x+y$ has size five, so
    $|C_{10}|=2\cdot5=10$. The allowed interval lengths are
    $1,2,5$. If a box $B$ has a full side, every parallel
    coordinate line contributes two points, giving $|C_{10}\cap B|=\frac{2|B|}{5}.$
    The remaining shapes have side lengths in $\{1,2\}$.
    A singleton has occupancy $0$ or $1$, an adjacent pair also
    has occupancy $0$ or $1$, since no two consecutive residues
    both belong to $S=\{0,2\}\subseteq\Z_5$.
    In a $2\times2$ box the coordinate sums are
    $a,a+1,a+1,a+2$, and its occupancy is
    \[
    n(a)=\mathbf1_S(a)+2\mathbf1_S(a+1)+\mathbf1_S(a+2),
    \qquad (n(0),n(1),n(2),n(3),n(4))=(2,2,1,1,2).
    \]
    Thus, according as $|B|=1,2,4$, the possible occupancies lie
    in $\{0,1\},\{0,1\},\{1,2\}$, respectively. Equivalently,
    \[
    |C_{10}\cap B|\in
       \{\lfloor2|B|/5\rfloor,\lceil2|B|/5\rceil\}
    \]
    for every allowed box. The complement inherits the same
    rounding condition at density $3/5$, and
    Theorem~\ref{thm:box} proves both assertions.
\end{proof}

    The restrictions on the character order matter, as the following
    five-point example shows.

    \begin{example}
        \label{ex:index-p-failure}
        In $\Z_5^2$, compare the coordinate sum fiber $C_5=\{(t,-t):t\in\Z_5\}$ with the line $K=\{(t,2t):t\in\Z_5\}$. Since both lines are subgroups, translation reduces their distance calculations
        to distances from the origin. On $C_5$ the four distances are $2\ell_5(t)$,
        $t=\pm1,\pm2$. On $K$, multiplication by two interchanges the
        two nonzero Lee lengths, so $\ell_5(t)+\ell_5(2t)=3$ for each
        $t\ne0$. Thus their enumerators are
        \[
        W_{C_5}(z)=1+2z^2+2z^4,\qquad W_K(z)=1+4z^3.
        \]
        Consequently, for every potential $f$ on $\{1,2,3,4\}$,
        \[
        E_f(C_5)-E_f(K)
        =2f(2)+2f(4)-4f(3)=2\Delta^2f(2).
        \]
        In particular, for the completely monotone potential
        $b_2(r)=\binom{4-r}{2}$,
        \[
        E_{b_2}(C_5)=2,\qquad E_{b_2}(K)=0.
        \]
        Thus even at prime index the coordinate sum fiber need not
        be universally optimal.
    \end{example}

    For odd side lengths and dimension at least two, Theorem~\ref{cor:index-three} and Proposition~\ref{prop:ten} provide the positive cases in the classification at the independence number. The next section completes this classification by treating the remaining parameters.

\section{Universally optimal sets in $\mathbb{Z}_q^n$ and maximum independent sets}
\label{sec:odd}

For the exponential potential $f_c(r)=c^r$, the
contribution of pairs at distance one dominates as $c\downarrow0$
whenever such pairs are present. Near $c=1$, the comparison is determined
by the first binomial moment at which two configurations differ. To
connect these two regimes, we use
the following classification of maximum independent sets due to Arras
and Joos~\cite[Lemmas~3.1, 3.2, and~3.4]{refAJ}.

\begin{lemma}\label{lem:stripes}
Let $q=2s+1\ge3$, $n\ge1$, and
$S=\{0,2,\ldots,2s-2\}\subset\Z_q$.
The distance one graph of $\Z_q^n$ has independence number $sq^{n-1}$.
Its maximum independent sets are precisely
\[
C_{a,\varepsilon}
=\biggl\{x\in\Z_q^n:
x_1+\sum_{i=2}^{n}\varepsilon_i x_i\in a+S\biggr\},
\qquad a\in\Z_q,\quad \varepsilon\in\{\pm1\}^{n-1}.
\]
Distinct parameter choices give distinct sets. Translations and
coordinate sign changes carry every such set to the stripe
$C=\{x:\sum_i x_i\in S\}$.
\end{lemma}

Thus all maximum independent sets have the same energy for any fixed
potential. For $q\ge5$ and $n\ge3$, we seek a competitor with lower
exponential energy near $c=1$ by changing one coefficient in the
defining sum of $C$ from $1$ to $2$. Both coefficients are units modulo the odd integer $q$, so
this change preserves the point counts in every proper coordinate
projection. These common projections force all binomial moments below
order $n$ to agree. The following lemma compares the first remaining
moment.

\begin{lemma}\label{lem:stripe-comparison}
Let $q=2s+1\ge5$, $n\ge3$, and $D=ns$. With $C$ and $S$ as above, put
\[
C'=\left\{x\in\Z_q^n:
x_1+\cdots+x_{n-1}+2x_n\in S\right\}.
\]
The sets $C$ and $C'$ both have size $M=sq^{n-1}$. Their $b_j$-energies
agree for $0\le j<n$, whereas $E_{b_n}(C')<E_{b_n}(C)$.
\end{lemma}

\begin{proof}
For $a\in\{1,2\}$, define
\[
\phi_a(x)=x_1+\cdots+x_{n-1}+ax_n,
\qquad A_a=\phi_a^{-1}(S),
\]
so $A_1=C$ and $A_2=C'$. We first compute the common projection
counts. The binomial expansion will then isolate a single product
kernel at order $n$. Its Fourier representation reduces the required
strict inequality to dimension three.

Let $J\subsetneq\{1,\ldots,n\}$, $m=|J|$, and
$\pi_J$ be the coordinate projection. For prescribed
$v\in\Z_q^J$ and $u\in S$, fix all but one of the coordinates
outside $J$. The coefficient of the remaining coordinate in $\phi_a$
is either $1$ or $2$, both of which are units modulo the odd integer
$q$. Hence
\begin{align*}
\#\{x:\pi_Jx=v,\ \phi_a(x)=u\}&=q^{n-m-1},\\
\#\{x\in A_a:\pi_Jx=v\}&=s q^{n-m-1}
\qquad(a=1,2).
\end{align*}
In particular, taking $J=\varnothing$ gives $|A_a|=M=sq^{n-1}$.
For every function $g:\Z_q^J\to\R$, these projection counts imply
\begin{align*}
\sum_{x,y\in A_a}g(\pi_Jx-\pi_Jy)
&=\sum_{v,w\in\Z_q^J}
  \#(A_a\cap\pi_J^{-1}(v))
  \#(A_a\cap\pi_J^{-1}(w))g(v-w)\\
&=s^2q^{2(n-m-1)}
  \sum_{v,w\in\Z_q^J}g(v-w).
\end{align*}
Both this expression and its diagonal contribution $Mg(0)$ are
independent of $a$. Thus the corresponding off-diagonal sums also
agree for $A_1$ and $A_2$.

Set $h(t)=s-\ell_q(t)$. The coefficient identity $(1+T)^{D-d(z)}=\prod_{i=1}^n(1+T)^{h(z_i)}$
gives
\[
b_j(d(z))=
\sum_{\substack{j_1+\cdots+j_n=j\\j_i\ge0}}
\prod_{i=1}^n\binom{h(z_i)}{j_i}.
\]
If $j<n$, every multi-index in the sum has a zero coordinate, so each
summand depends on a proper coordinate projection. The preceding
projection calculation therefore gives
\[
E_{b_j}(A_2)-E_{b_j}(A_1)=0\qquad(0\le j<n).
\]
For $j=n$, the only multi-index with all coordinates positive is
$(1,\ldots,1)$. After the other terms cancel, the energy difference
is therefore determined by the kernel $K(z)=\prod_{i=1}^n h(z_i).$

For $z\in\Z_q^n$, write
\[
N_a(z)=\#\{(x,y)\in A_a^2:x-y=z\}.
\]
The cancellation just proved yields
\[
\Delta_n:=E_{b_n}(C')-E_{b_n}(C)
=\frac1M\sum_{z\in\Z_q^n}K(z)\bigl(N_2(z)-N_1(z)\bigr).
\]
The term at $z=0$ vanishes because $N_1(0)=N_2(0)=M$.
To evaluate the remaining sum, let $\zeta=e^{2\pi i/q}$ and define,
with subscripts read modulo $q$,
\[
L_m=\left|\sum_{r=0}^{s-1}\zeta^{mr}\right|^2.
\]
If $u=\phi_a(y)$ and $v=\phi_a(x)$, then
$v-u=\phi_a(z)$. Each admissible pair $u,v\in S$ gives $q^{n-1}$
choices of $y$. Character orthogonality therefore gives
\begin{align*}
N_a(z)
&=q^{n-1}\sum_{u,v\in S}
  \mathbf1_{\{\phi_a(z)+u-v=0\}}\\
&=q^{n-2}\sum_{k=0}^{q-1}\zeta^{k\phi_a(z)}
  \left(\sum_{u\in S}\zeta^{ku}\right)
  \left(\sum_{v\in S}\zeta^{-kv}\right)\\
&=q^{n-2}\sum_{k=0}^{q-1}\zeta^{k\phi_a(z)}L_{2k}.
\end{align*}
Here the factor $L_{2k}$ follows from $S=\{2r:0\le r<s\}$.

For $I=\{0,\ldots,s-1\}$, the interval overlap formula gives
$h(t)=|I\cap(I+t)|$. Consequently,
\begin{align*}
\sum_{t\in\Z_q}h(t)\zeta^{mt}
=\sum_{r,r'\in I}\zeta^{m(r-r')}
=\left(\sum_{r\in I}\zeta^{mr}\right)
  \left(\sum_{r'\in I}\zeta^{-mr'}\right)=L_m.
\end{align*}
Substituting these formulas and separating the coordinate sums, with
$q^{n-2}/M=1/(qs)$, gives
\begin{align*}
\Delta_n
&=\frac1{qs}\sum_{k=0}^{q-1}L_{2k}
  \sum_{z\in\Z_q^n}K(z)
  \bigl(\zeta^{k\phi_2(z)}-\zeta^{k\phi_1(z)}\bigr)\\
&=\frac1{qs}\sum_{k=0}^{q-1}L_{2k}
  \left(\sum_t h(t)\zeta^{kt}\right)^{n-1}
  \left(\sum_t h(t)\zeta^{2kt}-\sum_t h(t)\zeta^{kt}\right).
\end{align*}
The $k=0$ term is zero. Thus
\begin{equation}\label{eq:5.4}
\Delta_n=\frac1{qs}\sum_{k=1}^{q-1}
L_{2k}L_k^{n-1}(L_{2k}-L_k).
\end{equation}

To prove that this sum is negative, we first show
$\Delta_n\le\Delta_3$. We then bound each summand for $n=3$ by a
quantity whose sum vanishes under the permutation $k\mapsto2k$.
Fix $1\le k<q$ and put
\[
\theta=\frac{\pi k}{q},\qquad x=L_k,\qquad y=L_{2k}.
\]
Because $2s=q-1$, we have $2s\theta=k\pi-\theta,
\
\cos(2s\theta)=(-1)^k\cos\theta,
\
\sin^2(2s\theta)=\sin^2\theta.$
The geometric-sum formula therefore becomes
\begin{align*}
x
&=\frac{\sin^2(s\theta)}{\sin^2\theta}
 =\frac{1-(-1)^k\cos\theta}{2(1-\cos^2\theta)}
 =\frac1{2(1+(-1)^k\cos\theta)},\\
y
&=\frac{\sin^2(2s\theta)}{\sin^2(2\theta)}
 =\frac{\sin^2\theta}{4\sin^2\theta\cos^2\theta}
 =\frac1{4\cos^2\theta}.
\end{align*}
All denominators are nonzero since $0<\theta<\pi$ and the oddness of $q$
implies $\theta\ne\pi/2$. In particular, $x>1/4$ and
$2x-1=-2x(-1)^k\cos\theta\ne0$. Eliminating $\cos\theta$ gives
\[
y=\frac{x^2}{(2x-1)^2},
\qquad
y-x=\frac{x(1-x)(4x-1)}{(2x-1)^2}.
\]
For $n\ge3$, we have $(x^{n-3}-1)(1-x)\le0.$
Multiplying this inequality by the positive factor
$x^3y(4x-1)/(2x-1)^2$ gives
\begin{align*}
x^{n-1}y(y-x)-x^2y(y-x)
&=\bigl(x^{n-3}-1\bigr)
  \frac{x^3y(1-x)(4x-1)}{(2x-1)^2}
\le0.
\end{align*}
Summing over $k$ therefore reduces the comparison to dimension three:
$\Delta_n\le\Delta_3$.

To bound $\Delta_3$, direct substitution yields
\[
2x^2y-x-y
=\frac{x(x-1)(2x^2-2x+1)}{(2x-1)^2}.
\]
Consequently,
\begin{align*}
x^2y(y-x)-\frac{y^2-x^2}{2}
=\frac12(2x^2y-x-y)(y-x)
=-\frac{x^2(x-1)^2(4x-1)(2x^2-2x+1)}
         {2(2x-1)^4}
\le0.
\end{align*}
The sign follows from $4x-1>0$ and
$2x^2-2x+1=2(x-1/2)^2+1/2>0$.
At $k=1$ the inequality is strict, because
\[
x=L_1=\frac1{4\sin^2(\pi/(2q))}>1
\qquad(q\ge5).
\]
Since doubling permutes $\Z_q\setminus\{0\}$, the comparison terms
cancel after summation. The strict inequality at $k=1$ then gives
\[
qs\,\Delta_3
<\frac12\sum_{k=1}^{q-1}(L_{2k}^2-L_k^2)=0.
\]
We conclude that $\Delta_n\le\Delta_3<0$ for every $n\ge3$.
\end{proof}

For odd $q\ge5$, the preceding lemma supplies a competitor when
$n\ge3$. In dimension two, a four-point exchange plays the same role,
except at $(q,n)=(5,2)$. Together, these comparisons give the following
classification.

\begin{theorem}\label{thm:odd}
Let $q=2s+1\ge5$, $n\ge2$, $M=sq^{n-1}$, and $f_c(r)=c^r$.
There exists an $M$-point set in $\Z_q^n$ that minimizes $E_{f_c}$ for
every $0<c<1$ if and only if $(q,n)=(5,2)$. In this exceptional case,
a universally optimal set exists.
\end{theorem}

\begin{proof}
Set $D=ns$. We first show that any common minimizer must be a maximum
independent set. Let $I$ be a maximum independent set. Then
\[
E_{f_c}(I)
=\frac1M\sum_{\substack{x,y\in I\\x\ne y}}c^{d_L(x,y)}
\le\frac{M(M-1)}M c^2=(M-1)c^2.
\]
If an $M$-point set $A$ contains an adjacent pair, the two corresponding
ordered pairs contribute
\[
E_{f_c}(A)\ge\frac{2c}{M}>(M-1)c^2
\qquad\left(0<c<\frac2{M(M-1)}\right).
\]
A common minimizer for all $c\in(0,1)$ must therefore be independent.
By Lemma~\ref{lem:stripes}, it is isometric to
$C=\{x:\sum_i x_i\in S\}$. To exclude a common minimizer outside the
exceptional case, it therefore suffices to find a set of the same
cardinality with lower energy than $C$ for some $c$.

To compare energies near $c=1$, put $t=c^{-1}-1>0$. For
$1\le r\le D$,
\begin{equation}\label{eq:exp-binomial}
c^r=c^D(1+t)^{D-r}
=c^D\sum_{j=0}^{D-1}b_j(r)t^j.
\end{equation}
For any $M$-point competitor $C'$, write
$\delta_j=E_{b_j}(C')-E_{b_j}(C)$. If
$\delta_0=\cdots=\delta_{m-1}=0$ and $\delta_m<0$, then
\begin{align*}
\frac{E_{f_c}(C')-E_{f_c}(C)}{c^Dt^m}
&=\delta_m+\sum_{j=m+1}^{D-1}\delta_jt^{j-m}<0
\qquad(\text{when }t\rightarrow0).
\end{align*}
Thus $E_{f_c}(C')<E_{f_c}(C)$ for all $c<1$ sufficiently close to one.
When $n\ge3$, Lemma~\ref{lem:stripe-comparison} supplies such a
competitor with $m=n$, proving nonexistence in these dimensions.

It remains to treat $n=2$ and $q\ge7$. Here we preserve the two
coordinate marginals by exchanging the corners of a rectangle.
Define
\[
p=(0,q-3),\qquad p'=(1,q-1),\qquad
w=(0,q-1),\qquad w'=(1,q-3),
\]
and set
\[
C'=(C\setminus\{p,p'\})\cup\{w,w'\},
\qquad
\eta=\mathbf1_{C'}-\mathbf1_C
=\mathbf1_{\{w\}}+\mathbf1_{\{w'\}}
 -\mathbf1_{\{p\}}-\mathbf1_{\{p'\}}.
\]
The coordinate sums show that $q-3,\ 0\in S,\ q-1,\ q-2\notin S,$
so the two removed points belong to $C$, the two added points do not,
and $|C'|=|C|=M$.

The exchange preserves both coordinate marginals. More explicitly,
\[
\sum_{v\in\Z_q}\eta(u,v)=0\quad(u\in\Z_q),
\qquad
\sum_{u\in\Z_q}\eta(u,v)=0\quad(v\in\Z_q).
\]
If $m_u=\#\{x\in C:x_1=u\}$, then for any $g:\Z_q\to\R$,
\[
\sum_{x,y\in C}g(x_1-y_1)
=\sum_{u,v\in\Z_q}m_um_vg(u-v)
=\sum_{x,y\in C'}g(x_1-y_1).
\]
The same identity holds for the second coordinate. Since the
cardinalities also agree, the diagonal contributions are equal, so
these identities remain valid for the off-diagonal pair sums.

For a displacement $(z_1,z_2)$ write
$u=\ell_q(z_1)$ and $v=\ell_q(z_2)$. The expansions
\begin{align*}
b_0(u+v)&=1,\\
b_1(u+v)&=2s-u-v,\\
b_2(u+v)
&=\binom{2s}{2}
 -\left(2s-\tfrac12\right)(u+v)
 +\tfrac12(u^2+v^2)+uv,
\end{align*}
therefore give $\delta_0=\delta_1=0$. All terms depending on a single
coordinate cancel, so only the mixed term $uv$ contributes to
$\delta_2$.

To evaluate this mixed contribution, put $\ell=\ell_q$ and define
\[
H(z_1,z_2)=\ell(z_1)\ell(z_2),\qquad
K(r)=\sum_{a\in\Z_q}\ell(a)\ell(r-a),\qquad
\Phi(r)=\sum_{t\in S}K(r-t).
\]
For $z=(z_1,z_2)$ with $z_1+z_2=r$, grouping the points of $C$ by
their coordinate sum gives
\begin{align*}
\sum_{x\in C}H(z-x)
=\sum_{t\in S}\sum_{x_1\in\Z_q}
  \ell(z_1-x_1)\ell(z_2-t+x_1)
=\sum_{t\in S}\sum_{a\in\Z_q}
  \ell(a)\ell(r-t-a)
=\Phi(r).
\end{align*}
Expanding $\mathbf1_{C'}=\mathbf1_C+\eta$ in the ordered pair sum,
and using $H(-z)=H(z)$ and $H(0)=0$, yields
\begin{align*}
M\delta_2
&=\sum_{z,z'}
\bigl(\mathbf1_{C'}(z)\mathbf1_{C'}(z')
     -\mathbf1_C(z)\mathbf1_C(z')\bigr)H(z-z')\\
&=2\sum_z\eta(z)\Phi(z_1+z_2)
 +\sum_{z,z'}\eta(z)\eta(z')H(z-z').
\end{align*}

Since $q-3-S=S$ and $K(-r)=K(r)$, we have
\begin{align*}
\Phi(q-3-r)
=\sum_{t\in S}K(q-3-r-t)
=\sum_{t'\in S}K(t'-r)=\Phi(r).
\end{align*}
In particular,
$\Phi(q-3)=\Phi(0)$ and $\Phi(q-2)=\Phi(q-1)$. Hence the linear term is
\begin{align*}
2\sum_z\eta(z)\Phi(z_1+z_2)
=2\bigl(\Phi(q-1)+\Phi(q-2)-\Phi(q-3)-\Phi(0)\bigr)
=4\bigl(\Phi(q-1)-\Phi(0)\bigr).
\end{align*}
The quadratic term depends only on the four exchanged points. In fact,
\[
H(p-p')=H(w-w')=2,
\qquad
H(a-b)=0
\quad(a\in\{p,p'\},\ b\in\{w,w'\}),
\]
because every mixed pair shares one coordinate. Thus
\[
\sum_{z,z'}\eta(z)\eta(z')H(z-z')
=2H(p-p')+2H(w-w')=8.
\]

To finish the comparison, we compute the difference
$\Phi(q-1)-\Phi(0)$. Evenness of $K$ gives
\begin{align*}
\Phi(q-1)-\Phi(0)
&=\sum_{i=0}^{s-1}\bigl(K(2i+1)-K(2i)\bigr)\\
&=K(1)-K(0)
 +\sum_{i=1}^{s-1}\bigl(K(2(s-i))-K(2i)\bigr)\\
&=K(1)-K(0).
\end{align*}
In the second line we used
$K(2i+1)=K(q-2i-1)=K(2(s-i))$.
Also,
\begin{align*}
K(0)-K(1)
&=\sum_a\ell(a)^2-\sum_a\ell(a)\ell(a-1)\\
&=\frac12\sum_a\bigl(\ell(a)-\ell(a-1)\bigr)^2\\
&=\frac12\bigl(s\cdot1^2+s\cdot(-1)^2+0^2\bigr)=s.
\end{align*}
The last equality follows from the cyclic distance list
$0,1,\ldots,s-1,s,s,s-1,\ldots,1$.
Combining the linear and quadratic contributions yields
\[
M\delta_2=4(K(1)-K(0))+8=8-4s=-2(q-5)<0.
\]
Thus $\delta_2$ is the first nonzero binomial-energy difference and is
negative. Equation~\eqref{eq:exp-binomial} therefore excludes a common
minimizer for the exponential family in the remaining negative cases.

For the remaining case $(q,n)=(5,2)$,
Proposition~\ref{prop:ten} supplies a universally optimal ten-point set.
It minimizes each exponential potential because $(-1)^j\Delta^j f_c(r)=c^r(1-c)^j>0
\ (0<c<1).$
\end{proof}

We now combine Theorem~\ref{thm:odd} with the constructions of
Section~\ref{sec:checkerboards} and the one dimensional case. This gives
the existence classification for all homogeneous tori and identifies
the minimizers for each fixed exponential potential.

\begin{corollary}\label{cor:existence}
Let $q\ge2$ and $n\ge1$, and put
$M=\lfloor q/2\rfloor q^{n-1}$. The following conditions are equivalent.
\begin{enumerate}[label=\textup{(\roman*)}]
\item There is a universally optimal $M$-point set.
\item There is an $M$-point set that minimizes $E_{f_c}$ for every
$0<c<1$, where $f_c(r)=c^r$.
\item Either $q$ is even, or $q=3$, or $n=1$, or $(q,n)=(5,2)$.
\end{enumerate}
Whenever these conditions hold, the universally optimal $M$-point sets
are exactly the maximum independent sets. They form a single Lee isometry
class. For every fixed $c\in(0,1)$, the minimizers of $E_{f_c}$ are
precisely the same sets.
In particular, there are exactly ten universally optimal ten-point
subsets of $\Z_5^2$.
\end{corollary}

\begin{proof}
Fix a coordinate direction. Each of its $q^{n-1}$ parallel lines
is a $q$-cycle when $q\ge3$, or a single edge when $q=2$. An
independent set $A$ therefore satisfies
\[
|A|=\sum_L|A\cap L|
\le q^{n-1}\lfloor q/2\rfloor=M.
\]
The checkerboards attain this bound for even $q$, and
Lemma~\ref{lem:stripes} supplies a stripe attaining it for odd $q$.

Every exponential $f_c$ is completely monotone, so
\textup{(i)} implies \textup{(ii)}. If \textup{(iii)} fails, then
\[
q\ge5\text{ is odd},\qquad n\ge2,\qquad(q,n)\ne(5,2),
\]
and Theorem~\ref{thm:odd} excludes \textup{(ii)}. Conversely,
Theorems~\ref{thm:checkerboard}, \ref{cor:index-three}, and
Proposition~\ref{prop:ten} prove \textup{(i)} when $q$ is even,
$q=3$, or $(q,n)=(5,2)$, respectively.

The only positive case not yet covered is $n=1$ with $q=2s+1\ge5$.
The stripe $S$ is a \emph{maximally even} set; its optimality for strictly
decreasing, strictly convex potentials also follows
from~\cite[Theorem~2.16]{BCL}. The following box count proves optimality
for every completely monotone sequence. The periodic lift of $S$ is
\[
\widetilde S=S+q\Z
=\left\{\left\lfloor\frac{jq}{s}\right\rfloor:j\in\Z\right\},
\]
since $j=bs+r$, $0\le r<s$, gives
$\lfloor jq/s\rfloor=bq+2r$. For integers $t$ and $1\le a\le q$,
\begin{align*}
\#\bigl(\widetilde S\cap\{t,\ldots,t+a-1\}\bigr)
&=\#\left\{j\in\Z:
\frac{ts}{q}\le j<\frac{(t+a)s}{q}\right\}\\
&=\left\lceil\frac{(t+a)s}{q}\right\rceil
 -\left\lceil\frac{ts}{q}\right\rceil\\
&\in\left\{\left\lfloor\frac{as}{q}\right\rfloor,
             \left\lceil\frac{as}{q}\right\rceil\right\}.
\end{align*}
Reduction modulo $q$ is injective on each such interval, so the same
count holds for every translated cyclic interval. Theorem~\ref{thm:box}
therefore applies at density $s/q$ and proves universal optimality of
$S$. This completes the equivalence of \textup{(i)}--\textup{(iii)}.

Having established existence, we now identify all universally optimal
sets in the positive cases. Put
$D=n\lfloor q/2\rfloor$. On the off-diagonal distances, $b_{D-1}(r)=\binom{D-r}{D-1}=\mathbf1_{\{r=1\}}.$
If $\nu_1(A)$ is the number of unordered adjacent pairs in $A$, then
\[
E_{b_{D-1}}(A)=\frac{2\nu_1(A)}M\ge0.
\]
An independent $M$-point set attains zero, so every universally optimal
set has $\nu_1(A)=0$ and is a maximum independent set. For odd $q$,
Lemma~\ref{lem:stripes} puts all maximum independent sets in one
Lee isometry class. For even $q$, the connected regular bipartite graph
argument in Theorem~\ref{thm:checkerboard} gives precisely the two
checkerboards, which differ by a translation. In each positive case, we
have already constructed a universally optimal set. Every maximum
independent set is isometric to this set and is therefore universally
optimal. The argument also includes $M=1$.
In $\Z_5^2$, the number of distinct parameter choices in
Lemma~\ref{lem:stripes} is $5\cdot2=10$.

It remains to identify the minimizers of each fixed exponential
potential. Fix a universally optimal set $U$ and $c\in(0,1)$. By
\eqref{eq:exp-binomial}, every $M$-point set $A$ satisfies
\begin{align*}
E_{f_c}(A)-E_{f_c}(U)
&=c^D\sum_{j=0}^{D-1}t^j
  \bigl(E_{b_j}(A)-E_{b_j}(U)\bigr),
\qquad t=c^{-1}-1>0.
\end{align*}
Universal optimality of $U$ makes every energy difference nonnegative,
and all weights $c^Dt^j$ are positive. Equality therefore holds exactly
when every fundamental energy agrees:
\begin{align*}
E_{f_c}(A)=E_{f_c}(U)
&\iff E_{b_j}(A)=E_{b_j}(U)\quad(0\le j<D)\\
&\iff A\text{ is universally optimal},
\end{align*}
where the last equivalence is Lemma~\ref{lem:fundamental}. This proves
the assertion for every fixed exponential potential.
\end{proof}

\section{Universally optimal $5^{n-1}$-point subsets in $\mathbb{Z}_5^n$}
    \label{sec:quinary}

    We now consider sets of density $1/5$ in $\Z_5^n$ and ask how to recognize universal optimality. We construct a certificate in the ordinary Fourier LP for each binomial potential, with a shared set of Fourier zeros, and exhibit a configuration attaining all the resulting bounds. We then show that the certificate for the potential supported at distances one and two has exactly this zero set. Its equality condition will therefore force equality for every binomial potential. The resulting support characterization also leads to constructions from Latin hypercubes.

    The certificates come from coefficient bounds for products of one coordinate generating polynomials. The two types of nonzero coordinate frequencies give two linear factors, and balancing their multiplicities will give the common zeros. A zero coordinate frequency contributes a convex combination of products of two such factors. This accounts for the range $N\ge n$ in the coefficient inequality below and allows the same bound to apply at every frequency.

    Fix $n\ge2$ and put
    \[
    \alpha=\frac{3+\sqrt5}{2},\qquad
    \beta=\frac{3-\sqrt5}{2},\qquad
    P(y)=1+3y+y^2=(1+\alpha y)(1+\beta y).
    \]
    Thus $\alpha+\beta=3$ and $\alpha\beta=1$. For $n=2k$ define
    \begin{equation*}
        A_m=[y^m]P(y)^k,\qquad
        c_m=[y^{m-1}](1+3y/2)P(y)^{k-1},
    \end{equation*}
    and for $n=2k+1$ define
    \begin{equation*}
        A_m=[y^m](1+3y/2)P(y)^k,\qquad
        c_m=[y^{m-1}]P(y)^k.
    \end{equation*}
    All these numbers are nonnegative, and $A_m=c_m=0$ for $m>n$.
    Coefficient extraction at a negative index is understood to give zero.

    \begin{lemma}
        \label{lem:spectral}
        For integers $N\ge n$, $0\le r\le N$, and $m\ge1$,
        \begin{equation}
            [y^m](1+\alpha y)^r(1+\beta y)^{N-r}
            \ge A_m+c_m\bigl(\alpha r+\beta(N-r)-3n/2\bigr).
            \label{eq:6.3}
        \end{equation}
        When $N=n$, equality holds at $r=k$ if $n=2k$, and at $r=k,k+1$
        if $n=2k+1$.
    \end{lemma}

    \begin{proof}
    If $m>n$, then $A_m=c_m=0$ and the left side is nonnegative.
    If $m=1$, then $A_1=3n/2$, $c_1=1$, and both sides equal
    $\alpha r+\beta(N-r)$. Hence assume $2\le m\le n$, and write
    \[
    U=1+\alpha y,\quad W=1+\beta y,\quad
    \delta=\alpha-\beta=\sqrt5,\quad T_0=\frac{3n}{2}.
    \]
    For $0\le r\le N$ put
    \[
    p_N(r)=[y^m]U^rW^{N-r},\qquad
    t_N(r)=\alpha r+\beta(N-r).
    \]
    Since $U-W=\delta y$, the slopes between consecutive points satisfy
    \begin{align*}
    D_{N,r}
      &:=\frac{p_N(r+1)-p_N(r)}{t_N(r+1)-t_N(r)}
        =[y^{m-1}]U^rW^{N-r-1},\\
    D_{N,r+1}-D_{N,r}
      &=\delta[y^{m-2}]U^rW^{N-r-2}\ge0
        \qquad(0\le r\le N-2).
    \end{align*}

    The slopes are therefore nondecreasing. We first use this fact to
    establish a supporting line when $N=n$. Put $k=\lfloor n/2\rfloor$.
    For $n=2k$, the identity $UW=P$ gives
    \begin{align*}
    t_n(k)&=T_0,& p_n(k)&=A_m,\\
    D_{n,k-1}&=[y^{m-1}]WP^{k-1},&
    D_{n,k}&=[y^{m-1}]UP^{k-1},\\
    \frac{D_{n,k-1}+D_{n,k}}2
      &=[y^{m-1}](1+3y/2)P^{k-1}=c_m.
    \end{align*}
    For $n=2k+1$, the corresponding identities are
    \begin{align*}
    t_n(k)&=T_0-\delta/2,& t_n(k+1)&=T_0+\delta/2,\\
    p_n(k)&=[y^m]WP^k=A_m-\delta c_m/2,\\
    p_n(k+1)&=[y^m]UP^k=A_m+\delta c_m/2,\\
    D_{n,k}&=[y^{m-1}]P^k=c_m.
    \end{align*}
    Together with monotonicity of the slopes, these identities give,
    in both parities,
    \[
    p_n(k)=A_m+c_m(t_n(k)-T_0),\qquad
    D_{n,i}\le c_m\ (i<k),\quad D_{n,i}\ge c_m\ (i\ge k).
    \]
    Telescoping expresses the difference from this line as
    \[
    p_n(r)-A_m-c_m(t_n(r)-T_0)
    =\begin{cases}
      \displaystyle\delta\sum_{i=k}^{r-1}(D_{n,i}-c_m),&r\ge k,\\[5pt]
      \displaystyle\delta\sum_{i=r}^{k-1}(c_m-D_{n,i}),&r<k,
    \end{cases}
    \]
    Each expression is nonnegative. The sums are empty at $r=k$;
    when $n$ is odd, the sum at $r=k+1$
    also vanishes because $D_{n,k}=c_m$.

    To extend the bound to $N>n$, we compare the minima after
    subtracting the linear term of slope $c_m$. Set
    \[
    F_N(r)=p_N(r)-c_m t_N(r),\qquad
    \mu_N=\min_{0\le r\le N}F_N(r).
    \]
    Coefficientwise $0\le W\le U$, and each expression defining $c_m$
    is a coefficient of a product of $n-1$ such factors, or the average
    of two such coefficients. Therefore
    \[
    c_m\le[y^{m-1}]U^{n-1}
        \le[y^{m-1}]U^{N-1}=D_{N,N-1}.
    \]
    Consequently $F_N(N)-F_N(N-1)=\delta(D_{N,N-1}-c_m)\ge0,$
    so a minimizing index $r$ can be chosen with $r<N$. At this index,
    \[
    0\le F_N(r+1)-F_N(r)=\delta(D_{N,r}-c_m).
    \]
    Since $r<N$, we may delete one factor $W=1+\beta y$, obtaining
    \begin{align*}
    p_N(r)-p_{N-1}(r)&=\beta D_{N,r},\\
    t_N(r)-t_{N-1}(r)&=\beta,\\
    \mu_N=F_N(r)
      &=F_{N-1}(r)+\beta(D_{N,r}-c_m)\ge\mu_{N-1}.
    \end{align*}
    Iterating this comparison and using equality at the central index
    for $N=n$ yields
    \[
    F_N(r)\ge\mu_N\ge\mu_n=A_m-T_0c_m,
    \]
    which is \eqref{eq:6.3}.
\end{proof}

    We next convert the coefficient bound into Fourier-LP certificates
    for the binomial potentials. For $u\in\Z_5$, define the
    one-coordinate generating polynomial
    \[
    Q_u(y)=\sum_{z\in\Z_5}y^{\ell_5(z)}(1+y)^{2-\ell_5(z)}\chi_u(z),
    \qquad a_u=2\cos(2\pi u/5).
    \]
    Grouping the terms by the Lee length of $z$ gives
    \[
    Q_u(y)=(1+y)^2+a_u y(1+y)+a_{2u}y^2
    =1+(2+a_u)y+(1+a_u+a_{2u})y^2.
    \]
    Here $1+a_u+a_{2u}=\sum_{z\in\Z_5}\chi_u(z)$ equals $5$ for $u=0$
    and $0$ otherwise. For $u\ne0$, $a_u^2=2+a_{2u}=1-a_u$.
    The signs $a_1=a_4>0$ and $a_2=a_3<0$ identify the two roots
    and give the following table.
    \begin{equation*}
        \begin{array}{c|c|c}
            u&Q_u(y)&a_u\\ \hline
            0&1+4y+5y^2&2\\
            1,4&1+\alpha y&\alpha-2\\
            2,3&1+\beta y&\beta-2
        \end{array}
    \end{equation*}
    To include zero coordinates in the coefficient bound, we express
    the zero-frequency polynomial as a convex combination:
    \begin{equation}
        1+4y+5y^2=w(1+\alpha y)^2+(1-w)(1+\beta y)^2,
        \qquad w=\frac{5+\sqrt5}{10}.
        \label{eq:6.5}
    \end{equation}
    The choice $w=(2-\beta)/(\alpha-\beta)$ makes the weighted mean of
    $\alpha,\beta$ equal to $2$. Since both satisfy $z^2=3z-1$, their
    weighted mean of the squares is $3\cdot2-1=5$, giving the quadratic
    coefficient in \eqref{eq:6.5}.

    For $u\in\Z_5^n$, put
    \[
    S(u)=\sum_i a_{u_i},\qquad
    r(u)=|\{i:u_i\in\{1,4\}\}|,\qquad
    s(u)=|\{i:u_i\in\{2,3\}\}|,
    \]
    and let $\ell=n_0(u)$ be the number of zero coordinates. Thus
    $r(u)+s(u)+\ell=n$. With $U=1+\alpha y$ and $W=1+\beta y$,
    \eqref{eq:6.5} gives the explicit convex combination
    \begin{align*}
    \prod_{i=1}^n Q_{u_i}(y)
      &=U^{r(u)}W^{s(u)}\bigl(wU^2+(1-w)W^2\bigr)^\ell\\
      &=\sum_{a=0}^\ell\omega_a
            U^{r(u)+2a}W^{s(u)+2(\ell-a)},
    \end{align*}
    where $\omega_a=\binom{\ell}{a}w^a(1-w)^{\ell-a}\ge0$ and $
    \sum_{a=0}^\ell\omega_a=1.$
    Each product in this combination has $N=n+\ell\ge n$ linear
    factors. The sum of their coefficients of $y$ is
    \[
    T_a=\alpha(r(u)+2a)+\beta(s(u)+2\ell-2a).
    \]
    Using $\sum_a a\omega_a=\ell w$ and
    $w\alpha+(1-w)\beta=2$, we obtain
    \begin{align*}
    \sum_a\omega_aT_a
       &=\alpha r(u)+\beta s(u)
                   +2\ell\bigl(w\alpha+(1-w)\beta\bigr)\\
       &=\alpha r(u)+\beta s(u)+4\ell
        =\sum_i(2+a_{u_i})=2n+S(u).
    \end{align*}
    Applying Lemma~\ref{lem:spectral} to each product and averaging gives
    \begin{align*}
    [y^m]\prod_iQ_{u_i}(y)
       &\ge\sum_a\omega_a\bigl(A_m+c_m(T_a-3n/2)\bigr)\\
       &=A_m+c_m\bigl(S(u)+n/2\bigr).
    \end{align*}

    We now identify these coefficients with Fourier transforms. Fix
    $0\le j\le2n-1$ and put $m=2n-j$. Since $d(x)=\sum_i\ell_5(x_i)$
    and $\chi_u(x)=\prod_i\chi_{u_i}(x_i)$, expanding the product gives
    \[
    \prod_{i=1}^n Q_{u_i}(y)
    =\sum_{x\in\Z_5^n} y^{d(x)}(1+y)^{2n-d(x)}\chi_u(x).
    \]
    In the term indexed by $x$, a power $y^{2n-j}$ is obtained by choosing
    the constant term in exactly $j$ of the $2n-d(x)$ linear factors. Thus
    \[
    [y^{2n-j}]\,y^{d(x)}(1+y)^{2n-d(x)}
    =\binom{2n-d(x)}{2n-j-d(x)}
    =\binom{2n-d(x)}{j}=b_j(d(x)).
    \]
    Extracting the coefficient term by term and using $d(-x)=d(x)$
    to substitute $x\mapsto -x$ in the character sum gives
    \[
    [y^m]\prod_i Q_{u_i}(y)
    =\sum_x b_j(d(x))\chi_u(x)
    =\sum_x b_j(d(x))\overline{\chi_u(x)}
    =5^n\widehat{b_j\circ d}(u).
    \]
    The distance-one displacements are $\pm e_i$, so their character
    sum is $\sum_i(\chi_u(e_i)+\chi_u(-e_i))=S(u)$. The character
    sum of the origin alone is $1$.
    Consequently, the auxiliary function
    \begin{equation*}
        h_j(x)=\binom{2n-d(x)}{j}
        -c_m\mathbf1_{\{d(x)=1\}}
        -(A_m+nc_m/2)\mathbf1_{\{x=0\}}
    \end{equation*}
    minorizes $b_j(d(x))$ for $x\ne0$ because $c_m\ge0$.
    Its Fourier transform satisfies
    \[
    5^n\widehat h_j(u)
    =[y^m]\prod_iQ_{u_i}(y)-A_m-c_m\bigl(S(u)+n/2\bigr)\ge0.
    \]
    The equality cases of Lemma~\ref{lem:spectral} provide common zeros
    of the transforms of these certificates. Define the set of balanced
    nonzero frequencies
    \[
    \Omega_n=\{u\in(\Z_5\setminus\{0\})^n:|r(u)-s(u)|\le1\}.
    \]
    At every $u\in\Omega_n$, the two type counts differ by at most one and
    $n_0(u)=0$. Lemma~\ref{lem:spectral} therefore gives $\widehat h_j(u)=0$
    for all $j$.

    We next relate the certificates to the first two distances.
    For a set $C$, let
    \[
    \nu_r(C)=|\{\{x,y\}\subseteq C:d_L(x,y)=r\}|,
    \qquad
    \Lambda(C)=(2n-1)\nu_1(C)+\nu_2(C).
    \]
    The potential $b_{2n-2}$ takes the values $2n-1$ and $1$ at distances
    one and two and vanishes at larger distances. Hence $E_{b_{2n-2}}(C)=\frac{2}{|C|}\Lambda(C).$

    The following theorem shows that this single potential detects
    universal optimality. Its equality condition is precisely the
    support restriction suggested by the common zeros above.

    \begin{theorem}
        \label{thm:quinary-characterization}
        Let $n\ge2$ and $C\subseteq\Z_5^n$ with $|C|=5^{n-1}$. The following conditions are equivalent.
        \begin{enumerate}[label=\textup{(\roman*)}]
            \item $C$ is universally optimal.
            \item $C$ is ordinary Fourier-LP universally optimal.
            \item $C$ minimizes $\Lambda(C)$ among all $5^{n-1}$-point subsets of $\Z_5^n$.
            \item The character sums of $C$ satisfy
            \begin{equation}
                S_C(u)=0\qquad
                \bigl(u\notin\{0\}\cup\Omega_n\bigr).
                \label{eq:quinary-support}
            \end{equation}
        \end{enumerate}
        Every set satisfying these conditions projects bijectively onto any $n-1$ coordinates and hence
        is an orthogonal array of strength $n-1$ and index one.
    \end{theorem}

    \begin{proof}
    We first show that all the Fourier bounds are attained, so that
    their equality conditions can be used to characterize the minimizers.
    Put $M=5^{n-1}$ and $V=5^n$. For $0\le j\le2n-1$, set
    $m=2n-j$ and
    \[
    B_j=M\widehat h_j(0)-h_j(0).
    \]
    Theorem~\ref{thm:fourier} bounds the ordinary Fourier-LP objective
    from below by $B_j$. For $x\ne0$, the definition of $h_j$ gives
    \[
    b_j(d(x))-h_j(x)=c_m\mathbf1_{\{d(x)=1\}}.
    \]
    Each unordered distance-one pair contributes twice to the first sum
    in \eqref{eq:fourier-gap}. Applying that identity with $f=b_j$ and
    $h=h_j$ therefore gives
    \begin{equation}
    E_{b_j}(C)-B_j
      =\frac{2c_m}{M}\nu_1(C)
       +\frac1M\sum_{u\ne0}\widehat h_j(u)|S_C(u)|^2\ge0.
    \label{eq:quinary-basic-gap}
    \end{equation}

    To attain all the bounds simultaneously, write
    $k=\lfloor n/2\rfloor$ and let
    \[
    v_*=
    \begin{cases}
       (1,2,\ldots,1,2),&n=2k,\\
       (1,2,\ldots,1,2,1),&n=2k+1,
    \end{cases}
    \qquad
    C_* =\{x\in\Z_5^n:\langle v_*,x\rangle=0\}.
    \]
    Every entry of $v_*$ is nonzero, so
    \[
    |C_*|=5^{n-1}=M,\qquad
    \langle v_*,\pm e_i\rangle=\pm(v_*)_i\ne0,
    \qquad \nu_1(C_*)=0.
    \]
    Character orthogonality determines its Fourier support:
    \begin{align*}
    S_{C_*}(u)
      =\frac15\sum_{t\in\Z_5}\sum_{x\in\Z_5^n}
                      \chi_{u+tv_*}(x)
      =\frac V5\sum_{t\in\Z_5}\mathbf1_{\{u+tv_*=0\}}
       =M\mathbf1_{\{u\in\langle v_*\rangle\}}.
    \end{align*}
    If $t\ne0$, every pair $(t,2t)$ contains one entry from
    $\{1,4\}$ and one from $\{2,3\}$; hence $tv_*\in\Omega_n$.
    Since $\widehat h_j=0$ on $\Omega_n$, both terms in
    \eqref{eq:quinary-basic-gap} vanish for $C_*$. Thus $C_*$ attains
    every bound $B_j$, both among sets and in the ordinary Fourier LP.

    Every actual difference distribution is LP feasible, and
    $E_{b_{2n-2}}(C)=2\Lambda(C)/M$. Therefore
    \[
    \textup{(ii)}\Longrightarrow\textup{(i)}
      \Longrightarrow\textup{(iii)}.
    \]
    To prove \textup{(iii)}$\Longrightarrow$\textup{(iv)}, we determine
    the complete zero set of the Fourier transform of the certificate
    with $m=2$. Strict positivity outside $\Omega_n$ will force the
    nontrivial character sums of a minimizing set to vanish there by
    \eqref{eq:quinary-basic-gap}, yielding the necessary support condition.
    For a fixed $u$, write
    \[
    z_i=[y]Q_{u_i}(y)=2+a_{u_i},\qquad
    T=\sum_i z_i=2n+S(u),\qquad T_0=3n/2.
    \]
    The identities $\alpha^2=3\alpha-1$ and $\beta^2=3\beta-1$
    imply
    \[
    z_i^2=3z_i-1+5\mathbf1_{\{u_i=0\}},\qquad
    \sum_i z_i^2=3T-n+5n_0(u).
    \]
    A quadratic term in the product comes either from a factor with
    zero coordinate frequency or from linear terms in two distinct factors.
    Hence
    \begin{align*}
    [y^2]\prod_iQ_{u_i}(y)
       =5n_0(u)+\sum_{i<j}z_i z_j
       =5n_0(u)+\frac12\left(T^2-\sum_i z_i^2\right)
       =g(T)+\frac52n_0(u), 
    \end{align*}
    where $g(T)=\frac12(T^2-3T+n).$
    Put $\varepsilon=\mathbf1_{\{n\text{ odd}\}}$. Direct coefficient
    extraction from the definitions of $A_m,c_m$ gives
    \begin{align*}
    A_2&=
       \begin{cases}
          k+9\binom{k}{2},&n=2k,\\
          k+9\binom{k}{2}+\frac92k,&n=2k+1,
       \end{cases}\\
    c_2&=\frac{3(n-1)}2=g'(T_0).
    \end{align*}
    Then $A_2=\frac{9n^2-14n+5\varepsilon}{8}
    =g(T_0)+\frac58\varepsilon.$
    Since $g(T)-g(T_0)-g'(T_0)(T-T_0)=\frac12(T-T_0)^2,$
    the Fourier transform of the certificate satisfies
    \begin{equation}
    5^n\widehat h_{2n-2}(u)
      =\frac12\left(\left(S(u)+\frac n2\right)^2
                  +5n_0(u)-\frac54\varepsilon\right).
    \label{eq:local-gap}
    \end{equation}
    If $n_0(u)\ge1$, this is at least $15/8>0$. If $n_0(u)=0$, then
    \[
    r(u)+s(u)=n,\qquad
    S(u)+\frac n2=\frac{\sqrt5}{2}(r(u)-s(u)),
    \]
    and consequently $5^n\widehat h_{2n-2}(u)
    =\frac58\bigl((r(u)-s(u))^2-\varepsilon\bigr).$
    As $r(u)-s(u)\equiv n\pmod2$, this vanishes precisely on
    $\Omega_n$. If $C$ minimizes $\Lambda$, it attains the same bound
    as $C_*$. Equation~\eqref{eq:quinary-basic-gap}, together with
    $c_2>0$, then yields
    \[
    \nu_1(C)=0,\qquad
    S_C(u)=0\quad(u\notin\{0\}\cup\Omega_n),
    \]
    proving \textup{(iv)}.

    Conversely, suppose \textup{(iv)} holds. We first apply the support
    condition to the coordinate projections. Fix a coordinate $i$ and
    let $\pi_i$ delete this coordinate. For $a\in\Z_5^{n-1}$ define
    \[
    N_i(a)=|\{x\in C:\pi_i(x)=a\}|.
    \]
    Let $\iota_i(v)$ insert a zero in coordinate $i$. Then
    \[
    \sum_aN_i(a)\chi_v(a)=S_C(\iota_i(v))
       =\begin{cases}M,&v=0,\\0,&v\ne0,\end{cases}
    \]
    because a nonzero frequency with a zero coordinate is outside
    $\Omega_n$. Fourier inversion therefore gives
    \[
    N_i(a)=\frac1{5^{n-1}}\sum_v
              S_C(\iota_i(v))\overline{\chi_v(a)}
           =\frac M{5^{n-1}}=1.
    \]
    Thus every such projection is bijective, proving the assertion
    about orthogonal arrays (cf.~\cite{XuWu}). In particular,
    \[
    x,y\in C,\quad x-y\in\Z_5e_i
       \quad\Longrightarrow\quad x=y,
    \]
    so $\nu_1(C)=0$. Condition \textup{(iv)} and the common zeros on
    $\Omega_n$ now make every term in \eqref{eq:quinary-basic-gap} zero.
    For a completely monotone potential $f=\sum_j\lambda_jb_j$ with
    $\lambda_j\ge0$, we obtain
    \[
    E_f(C)=\sum_j\lambda_jB_j
       \le\sum_j\lambda_j E_{b_j}(a)=E_f(a)
    \]
    for every feasible Fourier-LP vector $a$. This proves \textup{(ii)}
    and completes the equivalence.
\end{proof}

    The spectral criterion also applies when the box criterion fails,
    as the following example shows.

    \begin{example}\label{ex:box-not-necessary}
        Let $C=\{x\in\Z_5^3:x_1+2x_2=x_3\}$. Its nonzero Fourier support consists of the nonzero multiples of $(1,2,-1)$, all of which lie in $\Omega_3$. Thus Theorem~\ref{thm:quinary-characterization} gives universal optimality. The box $\{0,1\}\times\{0\}\times\{0,1\}$ contains both $(0,0,0)$ and $(1,0,1)$, whereas its mean occupancy is $4/5$. Its occupancy is therefore not one of the two integers bracketing the mean.
    \end{example}

    We now use the support characterization to construct further examples.
    Pairing coordinates will balance the two types of nonzero frequencies,
    while the Latin property of a base set will exclude frequencies with
    a zero coordinate.

    A Latin hypercube $Q\subseteq\Z_5^r$ contains exactly one point on every
    line obtained by fixing all but one coordinate. It has $5^{r-1}$ points and
    is the graph of an $(r-1)$-ary quasigroup when $r\ge2$. For $r=1$ it
    is a singleton. When $n=2k$, define
    \[
    \phi_n(x)=(x_1+2x_2,\ldots,x_{2k-1}+2x_{2k})\in\Z_5^k.
    \]
    For $n=2k+1$, append the unpaired coordinate $x_n$ to this vector, so that
    $\phi_n(x)\in\Z_5^{k+1}$. In each pair, prescribing its image and choosing the second coordinate determines the first. The unpaired coordinate, when present, is fixed by its image. Thus $\phi_n$ is onto and every fiber has $5^k$ points.

    \begin{theorem}
        \label{thm:latin}
        Let $n\ge2$. For every Latin hypercube $Q$ in the image of $\phi_n$, the
        lift $C_Q=\phi_n^{-1}(Q)$ has size $5^{n-1}$ and attains the ordinary
        Fourier LP for every completely monotone potential.
    \end{theorem}

    \begin{proof}
    Write $k=\lfloor n/2\rfloor$ and $r=\lceil n/2\rceil$.
    Every fiber of $\phi_n$ has $5^k$ points, so
    \[
    |C_Q|=5^k|Q|=5^k5^{r-1}=5^{n-1}.
    \]
    To verify the support condition, we first use the Latin property of $Q$.
    If $v\in\Z_5^r\setminus\{0\}$ has $v_j=0$, the projection deleting
    coordinate $j$ is bijective on $Q$. Therefore
    \[
    S_Q(v)=\sum_{a\in\Z_5^{r-1}}
                   \chi_{(v_1,\ldots,\widehat v_j,\ldots,v_r)}(a)=0.
    \]
    Here $\widehat v_j$ indicates the omitted coordinate. If $r=1$
    (equivalently $n=2$), $Q$ is a singleton and there is no nonzero
    $v$ with a zero coordinate, so this assertion is vacuous.

    Next we determine how the character sums change under the lift.
    The dual map is
    \[
    \phi_n^*v=
    \begin{cases}
      (v_1,2v_1,\ldots,v_k,2v_k),&n=2k,\\
      (v_1,2v_1,\ldots,v_k,2v_k,v_r),&n=2k+1.
    \end{cases}
    \]
    It satisfies $\chi_{\phi_n^*v}(x)=\chi_v(\phi_n(x))$.
    To identify the frequencies that survive summation over a fiber,
    observe that
    \[
    \ker\phi_n=
    \{(-2t_1,t_1,\ldots,-2t_k,t_k):t_i\in\Z_5\},
    \]
    with a final zero coordinate in odd dimension. Thus
    \begin{align*}
    \sum_{z\in\ker\phi_n}\chi_u(z)
       =\prod_{i=1}^k
          \sum_{t_i\in\Z_5}
          \exp\!\left(\frac{2\pi i}{5}
                (u_{2i}-2u_{2i-1})t_i\right)
       =5^k\prod_{i=1}^k
                   \mathbf1_{\{u_{2i}=2u_{2i-1}\}}.
    \end{align*}
    Summing each fiber separately consequently gives
    \[
    S_{C_Q}(u)=
    \begin{cases}
       5^kS_Q(v),&u=\phi_n^*v,\\
       0,&u\notin\operatorname{im}\phi_n^*.
    \end{cases}
    \]
    At a nonzero frequency with $S_{C_Q}(u)\ne0$, all coordinates of
    $v$ must be nonzero. Every pair $(v_i,2v_i)$ then contributes one
    coordinate of each type. Accordingly
    \[
    n_0(u)=0,\qquad
    (r(u),s(u))=
    \begin{cases}
       (k,k),&n=2k,\\
       (k+1,k)\text{ or }(k,k+1),&n=2k+1.
    \end{cases}
    \]
    Thus $u\in\Omega_n$, and Theorem~\ref{thm:quinary-characterization}
    proves ordinary Fourier-LP universal optimality of $C_Q$.
\end{proof}

    In the plane, the construction recovers a familiar equidistant set.
    Up to translation, the lift is $L=\{(t,2t):t\in\Z_5\}$, the quinary instance of Xu's equidistant Lee design~\cite[Proposition~1 and p.~1087]{ref24}. For $t\ne0$, the coordinates of $(t,2t)$ have Lee lengths $1$ and $2$ in some order, so all its off-diagonal distances are $3$.

    This equidistance extends optimality to all nonincreasing discretely
    convex potentials. Since $L$ minimizes the fundamental potential $b_1(d)=4-d$, it maximizes mean distance. More explicitly, if $\bar d(A)$ is the mean over the $20$ ordered off-diagonal pairs of a five-point set $A$, then $\bar d(A)\le3$. If $f$ is nonincreasing and discretely convex on $\{1,2,3,4\}$,
    its piecewise linear extension $\widetilde f$ is convex and nonincreasing.
    Jensen's inequality therefore gives
    \[
    \frac{E_f(A)}4
      =\frac1{20}\sum_{\substack{x,y\in A\\x\ne y}}
                         \widetilde f(d_L(x,y))
      \ge\widetilde f\bigl(\bar d(A)\bigr)
      \ge\widetilde f(3)=f(3)=\frac{E_f(L)}4.
    \]
    For exponential potentials, the same planar conclusion follows from Bilyk, Nagel, and Ruohoniemi~\cite[Corollary~2.5]{BNR}. After scaling by $1/5$ and reflecting the second coordinate, $L$ is the five-point Fibonacci lattice, and the product of the circle kernels $u\mapsto c^{5\min(u,1-u)}$ equals $c^{d_L}$ on the grid.

    In dimension four, the construction already yields a set that is
    not a subgroup coset.

    \begin{example}
        \label{ex:four-dimensional-noncoset}
        Let $n=4$, put $s=x_1+2x_2,\ t=x_3+2x_4,$
        and let $\pi$ be the permutation of $\Z_5$ with value list
        $(0,1,2,4,3)$. The set
        \[
        C_\pi=\{x\in\Z_5^4:t=\pi(s)\}
        \]
        is the Latin lift associated with the graph of $\pi$, so
        Theorem~\ref{thm:latin} shows that it has $125$ points and is universally optimal.
        Since $\pi(0)=0$, $\pi(1)=1$, and $\pi(2)=2$, the points
        $0$, $(1,0,1,0)$, and $(2,0,2,0)$ belong to $C_\pi$.
        The sum of the last two is $(3,0,3,0)$, which lies outside $C_\pi$
        because $\pi(3)=4$. A coset containing $0$ is a subgroup, so this
        failure of closure shows that $C_\pi$ is not a coset. Since the
        isometries in Lemma~\ref{lem:isometries} preserve the class of cosets
        in both directions, no Lee isometric image of $C_\pi$ is a coset.
    \end{example}

\section*{Acknowledgements}
    The author is partially supported by the National Natural Science Foundation of China (Grant No.~12571357) and the Natural Science Basic Research Program of Shaanxi (Program No.~2025JC-YBMS-048).

    The author used ChatGPT (GPT-6) in exploratory mathematical discussions and in revising the exposition. The author takes full responsibility for the manuscript.


\begin{thebibliography}{99}
        \small
        \setlength{\itemsep}{1.5pt plus .2pt}
        \bibitem{refAJ} P. Arras and F. Joos, Independent sets in discrete tori of odd sidelength, \emph{SIAM Journal on Discrete Mathematics} \textbf{39} (2025), no.~4, 1953--1973.
        \bibitem{ref1} A. Ashikhmin and A. Barg, Binomial moments of the distance distribution: bounds and applications, \emph{IEEE Transactions on Information Theory} \textbf{45} (1999), no.~2, 438--452.
        \bibitem{ref2} J. Astola, The Lee-scheme and bounds for Lee-codes, \emph{Cybernetics and Systems} \textbf{13} (1982), no.~4, 331--343.
        \bibitem{ref3} H. Astola and I. Tabus, On the linear programming bound for linear Lee codes, \emph{SpringerPlus} \textbf{5} (2016), article~246.
        \bibitem{ref4} A. Barg, Stolarsky's invariance principle for finite metric spaces, \emph{Mathematika} \textbf{67} (2021), no.~1, 158--186.
        \bibitem{BNR} D. Bilyk, N. Nagel, and I. Ruohoniemi, Minimizing point configurations for tensor product energies on the torus, \emph{Constructive Approximation} (2026), published online.
        \bibitem{ref5} N. Bouman, J. Draisma, and J. S. H. van Leeuwaarden, Energy minimization of repelling particles on a toric grid, \emph{SIAM Journal on Discrete Mathematics} \textbf{27} (2013), no.~3, 1295--1312.
        \bibitem{ref6} P. G. Boyvalenkov, P. D. Dragnev, D. P. Hardin, E. B. Saff, and M. M. Stoyanova, Energy bounds for codes and designs in Hamming spaces, \emph{Designs, Codes and Cryptography} \textbf{82} (2017), nos.~1--2, 411--433.
        \bibitem{ref7} D. Brosch and E. de Klerk, Minimum energy configurations on a toric lattice as a quadratic assignment problem, \emph{Discrete Optimization} \textbf{44} (2022), part~2, article~100612.
        \bibitem{BCL} N. Bushaw, B. Cody, and C. Leffler, Sets of vertices with extremal energy, \emph{Discrete Mathematics} \textbf{348} (2025), no.~8, article~114466.
        \bibitem{ref8} H. Cohn and N. Elkies, New upper bounds on sphere packings I, \emph{Annals of Mathematics} (2) \textbf{157} (2003), no.~2, 689--714.
        \bibitem{ref9} H. Cohn and A. Kumar, Universally optimal distribution of points on spheres, \emph{Journal of the American Mathematical Society} \textbf{20} (2007), no.~1, 99--148.
        \bibitem{ref10} H. Cohn, A. Kumar, S. D. Miller, D. Radchenko, and M. Viazovska, Universal optimality of the $E_8$ and Leech lattices and interpolation formulas, \emph{Annals of Mathematics} (2) \textbf{196} (2022), no.~3, 983--1082.
        \bibitem{ref12} H. Cohn and Y. Zhao, Energy-minimizing error-correcting codes, \emph{IEEE Transactions on Information Theory} \textbf{60} (2014), no.~12, 7442--7450.
        \bibitem{ref13} P. Delsarte, Bounds for unrestricted codes, by linear programming, \emph{Philips Research Reports} \textbf{27} (1972), 272--289.
        \bibitem{ref14} P. Delsarte, An algebraic approach to the association schemes of coding theory, \emph{Philips Research Reports Supplements} (1973), no.~10, vi+97~pp.
        \bibitem{GKP} R. L. Graham, D. E. Knuth, and O. Patashnik, \emph{Concrete Mathematics: A Foundation for Computer Science}, 2nd ed., Addison-Wesley, Reading, MA, 1994.
        \bibitem{HIK} R. Hammack, W. Imrich, and S. Klav\v{z}ar, \emph{Handbook of Product Graphs}, 2nd ed., CRC Press, Boca Raton, FL, 2011.
        \bibitem{ref17} A. R. Hammons, Jr., P. V. Kumar, A. R. Calderbank, N. J. A. Sloane, and P. Sol\'e, The $\Z_4$-linearity of Kerdock, Preparata, Goethals, and related codes, \emph{IEEE Transactions on Information Theory} \textbf{40} (1994), no.~2, 301--319.
        \bibitem{ref20} C. Y. Lee, Some properties of nonbinary error-correcting codes, \emph{IRE Transactions on Information Theory} \textbf{4} (1958), no.~2, 77--82.
        \bibitem{ref22} S. C. Polak, Semidefinite programming bounds for Lee codes, \emph{Discrete Mathematics} \textbf{342} (2019), no.~9, 2579--2589.
        \bibitem{ref23} P. Sol\'e, The Lee association scheme, in G. Cohen and P. Godlewski (eds.), \emph{Coding Theory and Applications}, Lecture Notes in Computer Science \textbf{311}, Springer, 1988, 45--55.
        \bibitem{ref24} H. Xu, Universally optimal designs for computer experiments, \emph{Statistica Sinica} \textbf{9} (1999), no.~4, 1083--1088.
        \bibitem{XuWu} H. Xu and C. F. J. Wu, Generalized minimum aberration for asymmetrical fractional factorial designs, \emph{Annals of Statistics} \textbf{29} (2001), no.~4, 1066--1077.
        \bibitem{ref25} V. A. Yudin, The minimum of potential energy of a system of point charges, \emph{Discrete Mathematics and Applications} \textbf{3} (1993), no.~1, 75--81.
    \end{thebibliography}
\end{document}